\documentclass[11pt]{article}
\usepackage{amsfonts,amssymb}
\usepackage{amsmath}
\usepackage{amsthm}
\usepackage[matrix,arrow,curve]{xy}
\usepackage{graphicx}
\usepackage[font=small,labelfont=bf]{caption}
\newcommand{\p}{\partial}
\newcommand{\e}{\varepsilon}

\newcommand{\R}{{\mathbb R}}
\newcommand{\A}{\mathcal A}
\newcommand{\U}{\mathcal U_T}
\newcommand{\Ui}{\mathcal U_\infty}

\newtheorem{theorem}{Theorem}
\newtheorem{lemma}{Lemma}
\newtheorem{prop}{Proposition}
\newtheorem{corollary}{Corollary}

\title{Optimal Control for Linear Systems with $L^1$-norm Cost}
\author{Andrei~Agrachev\thanks{SISSA, Trieste, agrachev@sissa.it} \and Bettina~Kazandjian\thanks{École Polytechnique, Sorbonne Université, Paris, bettina.kazandjian@polytechnique.edu}}
\date{}
\begin{document}
\maketitle

\begin{abstract}
We study $L^1$-optimal stabilization of linear systems with finite and infinite horizons. Main results concern the existence, uniqueness and structure of optimal solutions, and the robustness of optimal cost.
\end{abstract}

%\noindent Communicated by Laura Poggiolini.

\section{Introduction}

The $L^1$-norm of control as a cost attracted much less attention in the mathematical optimal control theory than the
$L^2$-norm or the traditional time-optimal problem. This cost has interesting peculiarities and is very relevant in common situations where optimal behavior should contain periods of the movement with a switched off control. An important motivation is the
aerospace navigation where this cost is proportional to the fuel consumption (see \cite{Ca,Te}).

In this paper, we consider the most simple setting: optimal stabilization of a linear system with control taking values in the Euclidean unit ball. We will see that the absolute value of optimal control is either 0 or 1 at every moment of time. This means that the cost is equal to the time spent with the activated control, while in the time-optimal problem the cost is the whole time of movement. We study both finite and infinite horizon problems. Main results concern the existence, uniqueness and structure of optimal solutions, and robustness of the optimal cost. The time-optimal control comes as a special case when
the horizon is fixed to the minimum time to reach the target; in this case, the control is always activated.

\section{Finite Horizon}

We consider optimal control problem for linear system
$$
\dot x(t)=Ax(t)+Bu(t),\quad x(t)\in\R^n,\ u(t)\in U\subset\R^m,\ 0\le t\le T. \eqno (1)
$$
Here $A:\R^n\to\R^n,\ B:\R^m\to\R^n,$ are linear maps, $U=\{u\in\R^m:|u|\le 1\}$ is a Euclidean unit ball, $|u|=\langle u,u\rangle^{\frac 12}$.
We would like to minimize the cost $J(u)=\int_0^T|u(t)|\,dt$ among all pairs $(u, x)$ which satisfy (1) and boundary conditions $x(0)=x_0,\ x(T)=0$. The time segment $T$ and the initial condition $x_0$ are parameters of the problem and
we are seeking for optimal solutions for all values of these parameters.

Let $E\subset\R^n$ be an invariant subspace of the operator $A$ and $Bv\in E,\ \forall\,v\in\R^m$.
Then $E$ is an invariant subspace of system (1) and the condition $x(T)=0$ implies that $x(t)\in E,\ 0\le t\le T.$
Minimal invariant subspace of $A$ which contains $\{Bv:v\in\R^m\}$ is the linear hull of $A^iBv,\ i=0,1,\ldots,n-1,\ v\in\R^m$. Indeed, $A^n$ is a linear combination of $A^i,\ 0\le i\le n-1$, according to the Cayley--Hamilton theorem.

We do not lose generality if we restrict our study to this invariant subspace or simply assume that this subspace is the whole $\R^n$. It is why
in this paper we always assume that $A, B$ satisfy the Kalman rank condition:
$$
span\{A^iBv: v\in\R^m,\ i=0,1,\ldots,n-1\}=\R^n. \eqno (2)
$$

Admissible controls $u(t),\ 0\le t\le T$ are measurable vector-functions with values in $U$. In other words, the space of admissible controls is the unit ball in $L^\infty([0,T];\R^m)$. We denote the space of admissible controls by $\mathcal U_T$ and we
equip $\mathcal U_T$ with the weak-$*$ topology: $u_n\rightharpoonup \bar u$ as $n\to\infty$ if and only if
$$
\int_0^T\langle u_n(t)-\bar u(t),v(t)\rangle\,dt\to 0,\quad\forall\,v\in L^1([0,T];\R^m).
$$
Then $\mathcal U_T$ is a compact topological space.

We denote by $x(t;u),\ 0\le t\le T,$ the solution of (1) with the initial condition $x(0)=x_0$. Cauchy formula for solutions of linear ordinary differential equations gives:
$$
x(t;u)=e^{tA}\left(x_0+\int_0^te^{-\tau A}Bu(\tau)\,d\tau\right).
$$
Hence $u\mapsto x(t;u)$ is a continuous linear map from $\mathcal U_T$ to $\R^n$.

Let
$$
\mathcal A_T=\{x_0\in\R^n: \exists\,u\in\mathcal U_T \mathrm{\ such\ that}\ x(T;u)=0\};
$$
then
$$
\mathcal A_T=\left\{\int_0^Te^{-tA}Bu(t)\,dt: u\in\mathcal U_T\right\}.
$$
It is easy to see that $\mathcal A_T,\ T>0,$ is a monotone growing family of compact convex subsets of $\R^n$, that are symmetric with respect to the origin.

\begin{prop} The origin of $\R^n$ is an interior point of $\mathcal A_T$ for any $T>0$.
\end{prop}
\noindent {\bf Proof.} Assume that $0\in\p\mathcal A_T$; the convexity of $\mathcal A_T$ implies that there exists
$\xi\in\R^n\setminus\{0\}$ such that $\langle\xi,x\rangle\le 0$ for any $x\in\mathcal A_T$. Hence
$\langle\xi,\alpha x\rangle\le 0$ for any $\alpha>0,\ x\in\mathcal A_T$. On the other hand, the rank condition (2) implies that
$$
\bigcup\limits_{\alpha\ge 1}\alpha\mathcal A_T=\left\{\int_0^Te^{-tA}Bu(t)\,dt: u\in\L^\infty([0,T];\R^m)\right\}=\R^n.
$$
This contradiction completes the proof. $\square$

\begin{corollary} $\A_T\subset int\A_{T+\varepsilon},\ \forall\,\varepsilon>0$.
\end{corollary}
\noindent Indeed, $\A_{T+\e}=\A_T+e^{-TA}\A_\e$.

\begin{prop} Let $x_0\in\mathcal A_T$; then there exists $\bar u\in\mathcal U_T$ such that
$$
J(\bar u)=\min\{J(u): u\in\mathcal U_T,\ x(T;u)=0\}.
$$
\end{prop}
\noindent{\bf Proof.} The desired result follows from the compactness of $\mathcal U_T$ and low semi-continuity of the functional $J : \mathcal U_T\to\R$ (for the topology introduced above). $\square $

The optimal control $\bar u$ may be not unique even for very simple systems. Indeed, let us consider the system
$\left\{\begin{aligned}\dot x^1&=x^2\\ \dot x^2&=u\end{aligned}\right.$, which describes a free particle on the line controlled by the external force.

Let $u\in\U,\ x(0;u)=x_0,\ x(T;u)=0$, where $x_0=(x_0^1,x_0^2)$. We have
$$
x_0^2=-\int_0^Tu(t)\,dt,\quad x_0^1=\int_0^Ttu(t)\,dt, \eqno (3)
$$
$$
J(u)=\int_0^T|u(t)|\,dt\ge \left|\int_0^Tu(t)\,dt\right|=|x_0^2|.
$$
Moreover, $J(u)=|x^2_0|$ if and only if control $u(t),\ 0\le t\le T,$ does not change sign. Hence any function from $\U$ which does not change sign is optimal for some initial condition; namely, for the initial condition (3).

Now assume that $\e\le u(t)\le 1-\e$, for some $\e>0$ and let $v\in L^\infty([0,T];\R)$ have zero average and zero first momentum: $$\int_0^Tv(t)\,dt=\int_0^Ttv(t)\,dt=0,$$ then $u+sv$ is an optimal control for our problem with the same initial condition for any $s$ sufficiently close to 0.

Coming back to the general case, we are going to characterize optimal controls by the Pontryagin maximum principle (see \cite{Po,AgSa}) and
to impose natural conditions on the pair $(A,B)$ which guarantee the uniqueness and simple structure of the optimal control.
Pontryagin maximum principle is a universal necessary optimality condition but it is also sufficient if the system is linear and if the constraints and the cost are convex, i.\,e. in our framework.

In order to formulate the maximum principle, we introduce the Hamiltonians of the system. Normal Hamiltonian:
$$h_u(p,x)=pAx+pBu-|u|,$$
where $p\in{\R^{n}}^*$ is a row, $p=(p_1,\ldots,p_n),\ px=p_1x^1+\cdots+p_nx^n$ (the product of a row and a column), and abnormal Hamiltonian: $h_u^0(p,x)= pAx+pBu$.

Both Hamiltonians depend on the parameter $u$. The Hamiltonian system for both Hamiltonians reads:
$$
\left\{\begin{aligned}\dot p&=-pA\\ \dot x&=Ax+Bu.\end{aligned}\right. 
$$
We say that $u\in\U$ is a normal extremal control if $x(T;u)=0$ and there exists a solution of the system $\dot p=-pAx$
such that
$$
h_{u(t)}(p(t),x(t;u))=\max\limits_{v\in U}h_v(p(t),x(t;u)),\quad \forall\,t\in[0,T].
$$
We say that $u$ is an abnormal extremal control if there exists a non-zero solution of the system $\dot p=-pAx$ such that
$$
h_{u(t)}^0(p(t),x(t;u))=\max\limits_{v\in U}h^0_v(p(t),x(t;u)),\quad \forall\,t\in[0,T].
$$

The Pontryagin maximum principle states that $u\in\U$, is optimal if and only if $x(T;u)=0$ and $u$ is a normal or abnormal extremal control.

Abnormal extremal controls do not depend on the cost $J$ and are actually extremal controls of the very well studied time-optimal problem. We briefly summarize their properties and devote the rest of the paper to the normal extremal controls.

\begin{theorem} Let $u\in\U,\ x(T;u)=0$. The control $u$ is an abnormal extremal control if and only if $x_0\notin\A_t,\ \forall\,t<T$. Moreover, any abnormal extremal control $u$ has the following properties:
\begin{enumerate}
\item $u$ is a unique solution of the equation $x(T;v)=0,\ v\in\U$;
\item $u$ is a piece-wise analytic vector-function with isolated jump discontinuities;
\item $|u(t)|=1,\ 0\le t\le T,$ and $u(t+0)=-u(t-0)$ in the discontinuity points.
\end{enumerate}
\end{theorem}

\noindent{\bf Proof.} We start from the proof of the statements {\it 2.} and {\it 3.}. The co-vector $p(t)$ from the maximality condition
has a form $p(t)=p^0e^{-tA}$, it is an analytic vector-function. Moreover, $u(t)^*=\frac 1{|p(t)B|}p(t)B$ when $p(t)B$ is not zero, where * is
the transposition, it transforms columns in rows and vice versa. The vector-function $p(t)B$ is analytic, it has only
isolated zeros if it is not identically zero. The identically zero case is excluded by the Kalman rank condition.

We have shown that $u$ is piece-wise analytic and takes values in the sphere $\partial U$. It remains to understand the structure of singularities. Let $p(\bar t)B=0$ and $k=\min\{i: \frac{d^i}{dt^i}p(\bar t)\ne 0\}$. We have:
$$
p(t)B=(t-\bar t)^ke+O\left((t-\bar t)^{k+1}\right).
$$
Then:
$$
u(t)=sign(t-\bar t)^k\frac e{|e|}+O(t-\bar t).
$$
Now we prove the first statement of the theorem. First of all, the Pontryagin maximum principle for the time-optimal problem implies that any time-optimal control must be an abnormal extremal control in our sense. Moreover, the equality
$x_0=-\int_0^Te^{-tA}Bu(t)\,dt$ implies:
$$
p(0)x_0=-\int_0^Tp(t)Bu(t)\,dt=-\int_0^T|p(t)B|\,dt<-\int_0^Tp(t)Bv(t)\,dt,
$$
for any $v\in\U$ which differs from $u$ on positive measure subset of $[0,T]$, and
$$
p(0)x_0<-\int_0^{t_1}p(t)Bv(t)\,dt, \quad \forall\,t_1<T, v\in\U.
$$

It follows that $u$ is a unique control which transfers $x_0$ to the origin in time $T$ and that $x_0\notin\A_{t_1},\ \forall\,t_1<T.\qquad \square$

\medskip
Now we analyze the maximality condition for the normal extremal controls. We have:
$$
h_{u(t)}(p(t),x)=\begin{cases} p(t)Ax+|p(t)B|-1&\text{if $|p(t)B|\ge 1$;}\\
p(t)Ax&\text{if $|p(t)B|\le 1$.}\end{cases}
$$

If $|p(t)B|>1$, then ${u(t)}^*=\frac 1{|p(t)B|}p(t)B$ like in the abnormal case; if $|p(t)B|<1$, then $u(t)=0$.
If $|p(t)B|=1$, then we can only say that ${u(t)}^*=\frac s{|p(t)B|}p(t)B$, for some $s\in[0,1]$.

Recall that $p(t)=p(0)e^{-tA}$; then either $|p(t)B|\equiv 1$ or the equation $|p(t)B|=1$ with unknown $t\in[0,T]$ has only
a finite number of solutions. If the first possibility is not realized, then optimal control is piece-wise analytic with jump discontinuities: it simply switches from $\frac 1{|p(t)B|}(p(t)B)^*$ to 0 and back when $|p(t)B|-1$ changes sign.

Moreover, optimal control is unique in this case as we see in the next proposition.

\begin{prop} Assume that for any $p^0\in{\R^n}^*$ there exists $t_0\in\R$ such that $|p^0e^{t_0A}|\ne 1$; then optimal control
is unique for any $T, x_0$.
\end{prop}

\noindent{\bf Proof.} Optimal control is either normal or abnormal extremal control. As we know (see Th.~1), abnormal extremal control is a unique control which transfers $x_0$ to the origin in time $T$. It corresponds to the pairs $T,x_0$
such that $x_0\in\p\A_T$. For all other pairs $T,x_0$ optimal controls are normal.

Assume that $u,w\in\U$ are two normal optimal controls for the same $T,x_0$. Then
$$
\int_0^Te^{-tA}Bu(t)\,dt=\int_0^Te^{-tA}Bw(t)\,dt,\quad \int_0^T|u(t)|\,dt=\int_0^T|w(t)|\,dt.
$$
Moreover, there exists $p^0\in{\R^n}^*$ such that
$$
p^0e^{-tA}Bu(t)-|u(t)|=\max\limits_{v\in U}(p^0e^{-tA}Bv-|v|),\quad \forall\,t\in[0,T], \eqno (4)
$$
and $u$ is a unique element of $\U$, which satisfies maximality condition (4).

In particular,
$$
p^0e^{-tA}Bu(t)-|u(t)|>p^0e^{-tA}Bw(t)-|w(t)|
$$
for any $t\in[0,T]$ such that $w(t)\ne u(t)$ and $|p^0e^{-tA}B|\ne 1$.

It follows that $u(t)=w(t)$ almost everywhere.\qquad$\square$

The assumption of Proposition~3 is violated for the considered above system $\left\{\begin{aligned}\dot x^1&=x^2\\ \dot x^2&=u\end{aligned}\right.$, where, as we know, optimal control is not unique.

\begin{lemma} Let $m=1$. The assumption of Proposition~3 is valid if and only if $\det (A)\ne 0$.
\end{lemma}

\noindent{\bf Proof.} We have $B:\R\to\R^n,\ B1=b\in\R^n$. In this case, the assumption of Proposition~3 can be rewritten
as follows: $p^0e^{tA}b\not\equiv const$ for any $p^0\in{\R^n}\setminus\{0\}$ or, equivalently,
$p^0Ae^{tA}b\not\equiv 0$. In other words, this assumption is valid if and only if
$A(span\{e^{tA}b:t\in\R\})=\R^n$.

The Kalman rank condition guaranties that $span\{e^{tA}b:t\in\R\}=\R^n$. Hence the assumption of Proposition~3 is valid if
and only if the operator $A$ verifies $\det(A)\neq0$.

Unfortunately, the last test does not work for $m\ge 2$. Here is a simple counter-example for $n=m=2$. We identify $\R^2$
with the complex plane ${\mathbb C}$. Complex numbers are identified with linear operators on ${\mathbb C}=\R^2$ acting
by the complex multiplication. We set $A=i,\ B=1$; then $e^{tA}B=e^{ti}$. If $p^0=(1,1)\in{\R^2}^*$, then
$|p^0e^{tA}B|=|e^{ti}|=1$. $\square$

To make things work for any $m$, we require the hyperbolicity of the operator $A$. Recall that $A$ is called hyperbolic if any eigenvalue of $A$ has a non-zero real part. A key property of the hyperbolic operators is as follows:
if $A$ is hyperbolic and $b\in\R^n\setminus\{0\}$, then either $|e^{tA}b|\to\infty$ as $t\to+\infty$ or
$|e^{tA}b|\to\infty$ as $t\to-\infty$. In the first case, $|e^{tA}b|\to 0$ as $t\to-\infty$ and in the second case,
$|e^{tA}b|\to 0$ as $t\to+\infty$.

\begin{lemma} If $A$ is hyperbolic, then the assumption of Proposition~3 is valid and optimal controls are unique.
\end{lemma}

\noindent{\bf Proof.} Let $e_1,\ldots,e_m$ be an orthonormal basis of $\R^m$, $b_i=Be_i,\ i=1,\ldots,m$. We have:
$|p^0e^{tA}B|^2=\sum\limits_{i=1}^m|p^0e^{tA}b_i|^2$. The key property of the hyperbolic operators implies that
$|p^0e^{tA}B|^2=0$ if and only if $p^0e^{tA}b_i=0,\ i=1,\ldots,m$. The differentiation of the last identities gives:
$p^0A^jb_i=0, \forall\,i,j$, but this is not possible for $p^0\ne 0$ due to the Kalman rank condition.
\qquad$\square$

Let $x_0\in\A_T$ and
$$
\mu_T(x_0)=\min\{J(u):u\in\U,\ x(T;u)=0\}
$$
be the optimal cost. In what follows, we are interested not only in the dependence of $\mu_T(x_0)$ on $T$ and $x_0$ but
also on the matrices $A,B$ and we use notations $\A_T=\A_T(A,B),\ \mu_T(x_0)=\mu_T(x_0;A,B),\
(T,x_0;A,B)\in\R_+\times\R^n\times\R^{n^2}\times\R^{nm}$.

\begin{theorem} The set
$$
\{(T,x_0;A,B): x_0\in int\A_T(A,B)\} \eqno (5)
$$
is an open subset of $\R_+\times\R^n\times\R^{n^2}\times\R^{nm}$ and the function
$$
(T,x_0;A,B)\mapsto\mu_T(x_0;A,B)
$$
is continuous on (5).
\end{theorem}

\noindent{\bf Proof.} In what follows, we assume that $\U\subset {\mathcal U}_{T'}$ for any $T'>T$, where
$u\in L^\infty([0,T];\R^m)$ is extended to the interval $(T,T']$ by zero.

First we prove the openness of set (5). Let $x_0\in int\A_T(A,B)$; take a simplex with vertices
$y_0,\ldots,y_n\in\A_T(A,B)$ such that $x_0$ is an interior point of this simplex. We have:
$$
y_i=\int_0^Te^{-tA}Bu_i\,dt,\quad u_i\in\U,\ i=0,\ldots,n.
$$
If $T',A',B'$ are close to $T,A,B$, then $y'_i=\int_0^{T'}e^{-tA'}B'u_i\,dt$ are close to $y_i,\ i=0,\ldots,n$, and $x_0$
is an interior point of the simplex with vertices $y'_i$, i.\,e. $x_0=\sum\limits_{i=0}^n\alpha_iy'_i$, where $\alpha_i\in(0,1)$. Hence
$$
x_0=\int_0^{T'}e^{-tA'}B'\sum\limits_{i=0}^n\alpha_iu_i(t)\,dt
$$
and $x_0\in int\A_{T'}(A',B')$.

The next step is the continuity of $\mu_T(x_0;A,B)$ with respect to $T$. First of all, $T\mapsto\mu_T(x_0)$ is a monotone decreasing function. Indeed, let $u\in\U$ be the minimizing control, $\mu_T(x_0)=J_T(u)$; then
$$
\mu_{T'}(x_0)\le J_{T'}(u)=J_{T}(u)=\mu_T(x_0),\quad \forall\,T'\ge T.
$$
The monotonicity implies the existence of the right and left limits:
$$
\lim\limits_{T'\searrow T}\mu_{T'}(x_0)=\mu_{T+0}(x_0),\quad \lim\limits_{T'\nearrow T}\mu_{T'}(x_0)=\mu_{T-0}(x_0).
$$

It is easy to see that $\mu_{T+0}(x_0)=\mu_T(x_0)$. Indeed, let $T_n\searrow T$ and $\mu_{T_n}(x_0)=J_{T_n}(u_n)$.
We may assume that $u_n\rightharpoonup\bar u$ as $n\to\infty$ by taking a sub-sequence if necessary; then
$\lim\limits_{n\to\infty}\mu_{T_n}(x_0)=J_T(\bar u)\ge\mu_T(x_0)$. On the other hand, $\mu_{T+0}(x_0)\le\mu_T(x_0)$;
hence $\mu_{T+0}(x_0)=\mu_T(x_0)$.

It remains to prove that $\mu_{T-0}(x_0)=\mu_T(x_0)$. This is more complicated. Let $\tilde u\in\U$ be the optimal control:
$J(\tilde u)=\mu_T(x_0)$ and $-\int_0^Te^{-tA}B\tilde u(t)\,dt=x_0$. We have to show that there exists an arbitrarily close to $\tilde u$ in the $L^1$-norm control $u'\in\U$ such that $-\int_0^{T'}e^{-tA}B u'(t)\,dt=x_0$ for some $T'<T$.

Recall that $\A_T$ is a full-dimensional convex compact which is symmetric with respect to the origin. Moreover, $x_0\in int\A_T$; hence $x_0\in\nu_0\A_T$ for some $\nu_0\in(0,1)$. Let $v_0\in\U$ be such that
$-\nu_0\int_0^Te^{-tA}Bv_0(t)\,dt=x_0$. Then
$$
-\int_0^Te^{-tA}B(s\nu_0v_0(t)+(1-s)\tilde u(t))\,dt=x_0,\quad \forall\,s\in[0,1].
$$
Let $v=s\nu_0v_0(t)+(1-s)\tilde u(t)$, where $s>0$ is so close to 1 that the $L^1$-norm of $v-\tilde u$ is smaller than a
preliminary chosen $\e>0$.

By construction, $v\in\nu\U$ for some $\nu\in(0,1)$ and $-\int_0^Te^{-tA}Bv(t)\,dt=x_0$.

Let $\tau\in[0,T)$, we set $\A_{\tau,T}=\left\{\int_\tau^Te^{-tA}Bu(t)\,dt: u\in\U\right\}$; then $\A_{\tau,T}$ is a
full-dimensional convex compact which is symmetric with respect to the origin. We have:
$$
x_0+\int_0^\tau e^{-tA}Bv(t)\,dt=-\int_\tau^Te^{-tA}Bv(t)\,dt\in\nu\A_{\tau,T}.
$$
Hence $-\int_\tau^Te^{-tA}Bv(t)\,dt\in int\A_{\tau,T}$. We claim that any interior point of $\A_{\tau,T}$ belongs to
$\A_{\tau,T_\tau}$ for some $T_\tau<T$. This fact can be proved in the same way as the openness of set (5), we leave the details to the reader.

Now we know that $x_0+\int_0^\tau e^{-tA}Bv(t)\,dt\in\A_{\tau,T_\tau}$ for some $T_\tau\in[\tau,T)$. Hence
$$
x_0+\int_0^\tau e^{-tA}Bv(t)\,dt=-\int_\tau^{T_\tau} e^{-tA}Bv_\tau(t)\,dt,
$$
for some $v_\tau\in\U$. We set
$$
u_\tau(t)=\begin{cases}v(t),&\text{if $0\le t\le\tau$;}\\ v_\tau(t),&\text{if $\tau<t\le T$.}\end{cases}
$$
Then $-\int_0^{T_\tau} e^{-tA}Bu_\tau(t)\,dt=x_0$ and $u_\tau$ tends to $v$ in the $L^1$-norm as $\tau\to T$.
These relations complete the proof of the continuity of $\mu_T(x_0;A,B)$ with respect to $T$.

The continuity with respect to $T$ will help us to prove the continuity with respect to all variables.
Let $(T_n,x_n;A_n,B_n)\to(T,x_0;A,B)$ as $n\to\infty$. We will separately prove the inequalities:
$$
\mu_T(x_0;A,B)\le\liminf\limits_{n\to\infty}\mu_{T_n}(x_n;A_n,B_n), \quad \mu_T(x_0;A,B)\ge\limsup\limits_{n\to\infty}\mu_{T_n}(x_n;A_n,B_n).
$$

The first inequality is easy. Let
$$
J_{T_n}(u_n)=\mu_{T_n}(x_n;A_n,B_n), \quad -\int_0^{T_n} e^{-tA_n}B_nu_n(t)\,dt=x_n.
$$
Let $u_{n_k}$ be a convergent sub-sequence, $u_{n_k}\rightharpoonup \bar u$ as $k\to\infty$, then
$$
J_T(\bar u)=\lim\limits_{k\to\infty}\mu_{T_{n_k}}(x_{n_k};A_{n_k},B_{n_k}),\quad -\int_0^T e^{-tA}B\bar u(t)\,dt=x_0
$$
and $\mu_T(x_0;A,B)\le J_T(\bar u)$. Hence $\mu_T(x_0;A,B)\le\liminf\limits_{n\to\infty}\mu_{T_n}(x_n;A_n,B_n)$.

Let us prove the second inequality. Take a small $\delta>0$; let
$$
J_{T-\delta}(u_\delta)=\mu_{T-\delta}(x_0;A,B), \quad -\int_0^{T-\delta} e^{-tA}Bu_\delta(t)\,dt=x_0.
$$
We set $y_n=x_n+\int_0^{T-\delta} e^{-tA_n}B_nu_\delta(t)\,dt$; then $y_n \rightarrow 0$. We claim that there exists a neighborhood of $0$ contained in $\underset{n\geq N}{\bigcap}\A_{T-\delta,T_n}(A_n,B_n)$ for some integer $N$. This fact can be proved in the same way as the openness of set (5), we leave the details to the reader. We consider $n$ large enough such that $|T_n-T|<\delta$ and $y_n\in\A_{T-\delta,T_n}(A_n,B_n)$. In particular $y_n=-\int_{T-\delta}^{T_n} e^{-tA_n}B_nv_n(t)\,dt$ for some $v_n\in{\mathcal U}_{T+\delta}$.
Let
$$
\hat u_n(t)=\begin{cases}u_\delta(t),&\text{if $0\le t\le T-\delta$;}\\ v_n(t),&\text{if $T-\delta<t\le T_n$.}\end{cases}
$$
Then $-\int_0^{T_n}e^{-tA_n}B_n\hat u_n(t)\,dt=x_n$ and
$$
\mu_{T_n}(x_n;A_n,B_n)\le J_{T_n}(\hat u_n)\le\mu_{T-\delta}(x_0;A,B)+2\delta.
$$
Hence $\limsup\limits_{n\to\infty}\mu_{T_n}(x_n;A_n,B_n)\le\mu_{T-\delta}(x_0;A,B)+2\delta$, for any $\delta>0$. It remains
to go to the limit as $\delta\to 0.\qquad\square$

\section{Infinite Horizon}

Let
$$
\Ui=\{u\in L^1([0,\infty);\R^m): |u(t)|\le 1,\ \forall\,t\ge 0\},\quad J_\infty(u)=\int_0^\infty|u(t)|\,dt,
$$
$\A_\infty=\bigcup\limits_{T>0}\A_T$; then $\A_\infty$ is an open convex subset of $\R^n$. Given $x_0\in\A_\infty$, we set
$$
\mu_\infty(x_0)=\inf\left\{J_\infty(u): u\in\Ui,\ \lim\limits_{t\to\infty}x(t;u)=0\right\}.
$$

\begin{prop} $\mu_\infty(x_0)=\inf\limits_{T>0}\mu_T(x_0)$.
\end{prop}

\noindent{\bf Proof.} We have: $\mu_T(x_0)=\inf\{J_T(u): u\in\U,\ x(T;u)=0\}$. The inclusion $\U\subset\Ui$ implies that
$\mu_T(x_0)\ge\mu_\infty(x_0),\forall\, T>0$. On the other hand, for any $\varepsilon>0$ and $u\in\Ui$, the relation
$x(t;u)\to0\ (t\to\infty)$ implies the existence of $T>0$ such that $x(T-\varepsilon;u)\in\A_{T-\e,T}$. It follows that
$x(T-\e;u)=-\int\limits_{T-\e}^Te^{-tA}Bv(t)\,dt$ for some $v\in\U$. We set:
 $$
u_\e(t)=\begin{cases}u(t)&\text{if\ $0\le t\le T-\e$;}\\ v(t)&\text{if\ $T-\e<t\le T$.}\end{cases}
$$
Then $x(T;u_\e)=0,\ J_T(u_\e)\le J_\infty(u)+\e$. Hence $\mu_T(x_0)\le\mu_\infty(x_0)+\e.\qquad\square$

In general, we cannot substitute the $\inf$ by the $\min$ in the definition of $\mu_\infty$ (see the next section for a simple counter-example). We can do it if $A$ is a hyperbolic operator. Moreover, the infinite horizon problem is reduced to
a finite horizon one in this case.

Let $A:\R^n\to\R^n$ be hyperbolic, then $\R^n=E^+\oplus E^-$, where $AE^\pm\subset E^\pm$ and invariant subspace
$E^+\ (E^-)$ corresponds to the eigenvalues of $A$ with positive (negative) real parts. We have $e^{tA}E^\pm=E^\pm$;
moreover, there exists $\alpha>0$ such that
$$
|e^{tA}x^+|\ge c_+e^{t\alpha}|x^+|,\quad  |e^{tA}x^-|\le c_-e^{-t\alpha}|x^-|,
$$
for some $c_\pm>0$ and any $x^\pm\in E^\pm,\ t\ge 0$.

Given $x\in\R^n$, we set $x=x^++x^-$, where $x^\pm\in E^\pm$ and define a linear map $B^+:\R^m\to E^+$ by the formula
$B^+u=(Bu)^+$. Finally, we set $\mu^+_T(x^+_0)=\mu_T(x^+_0;A,B^+)$ and $\mu^+_\infty(x_0)=\mu_\infty(x^+_0;A,B^+)$, the optimal cost for the system
$\dot x^+=Ax^++B^+u$ on $E^+$.

\begin{theorem} Let $A$ be hyperbolic; then for any $x_0\in\A_\infty$ there exists a unique optimal control $u\in\Ui$
for the infinite horizon problem with the initial condition $x_0$. Moreover, there exists $T>0$ such that
$\mu_\infty(x_0)=\mu^+_T(x^+_0)=J_T(u)$ and $u(t)=0$ for any $t>T$.
\end{theorem}

\noindent{\bf Proof.} It is easy to see that $t\mapsto x(t;u)^+$ is a solution of the system $\dot x^+=Ax^++B^+u$ on $E^+$,
for any $u\in\Ui$. Hence $\mu_\infty^+(x^+_0)\le\mu_\infty(x_0)$. Moreover, if $x(T;u)^+=0$, then $x(T;u)$ belongs to the asymptotically stable subspace of the system $\dot x=Ax$ and zero control transfers $x(T;u)$ to the origin in infinite time without augmentation of the cost.
Hence $\mu_\infty^+(x^+_0)\le\mu_\infty(x_0)\le\mu_T^+(x^+_0)$.

It remains to show that for any $x_0^+\in E^+$ there exists $T>0$ such that $\mu_\infty^+(x^+_0)=\mu_T^+(x^+_0)$.
To do that, we may assume that $E^+=\R^n$ in order to simplify notations. We fix $x_0\in\A_{T_0}$ and take $u_T\in\U$
such that $J_T(u_T)=\mu_T(x_0)$, for any $T>T_0$. Then $u_T$ is a normal extremal control and there exists $p_T\in\R^{n*}$ such that:
$$
|u_T(t)|=\begin{cases}1&\text{if\ $|p_Te^{-tA}B|>1$;}\\ 0&\text{if\ $|p_Te^{-tA}B|<1$,}\end{cases} \qquad0\le t\le T.
$$
Recall that $|p_Te^{-tA}B|\le ce^{-\alpha t}|p_T|$, for some positive constants $\alpha$ and $c$.

If $|p_T|$ are uniformly bounded, i.\,e. $|p_T|\le c',\ \forall\,T>T_0$, then $u_T(t)=0$ for any
$t>\frac 1\alpha(\ln c + \ln c')$ and we obtain that
$$
\mu_\infty(x_0)=\mu_T(x_0), \quad \forall\,T>\frac 1\alpha(\ln c + \ln c').
$$

If $|p_T|$ is not uniformly bounded, then there exists a sequence $T_k\to\infty,\ (k\to\infty)$ such that $|p_{T_k}|\to\infty$;
moreover, we may assume that $\frac 1{|p_{T_k}|}p_{T_k}\to\xi$, where $|\xi|=1$. Let us show that this is not possible.

We consider the set $\mathcal T=\{t\in[0,2T_0]: \xi e^{-tA}B=0\}$; this is a finite subset of $[0,2T_0]$. Let $O_\e\mathcal T$ be the radius $\e$ neighborhood of $\mathcal T$, this is the union of $\#\mathcal T$ intervals where each interval has length $2\e$. We fix a small enough $\e$ to guarantee that the measure of $O_\e\mathcal T$ is smaller than $T_0$.

There exists $\delta>0$ such that $|p_{T_k}e^{-tA}B|\ge\delta|p_{T_k}|$ for any $t\in[0,2T_0]\setminus O_\e\mathcal T$ and
any $k>\frac 1\delta$. Hence $|p_{T_k}e^{-tA}B|\ge 1$ for all $t\in[0,2T_0]\setminus O_\e\mathcal T$ if $k$ is sufficiently large. It follows that $|u_{T_k}(t)|=1,\ \forall\,t\in[0,2T_0]\setminus O_\e\mathcal T$ and
$\mu_{T_k}(x_0)=J(u_{T_k})>T_0$.

On the other hand, $\mu_{T_k}(x_0)\le\mu_{T_0}(x_0)\le T_0$. This contradiction proves that $|p_T|$ is uniformly bounded and
thus completes the proof of the theorem.

Suppose that $u_1,u_2 \in \Ui$ are optimal. Then there exists $T_1,T_2 > 0$ such that $\mu_\infty(x_0) = J_{T_i}({u_i}_{|[0;T_i]}) $ and $u_i(t)=0, \forall t\geq T_i$. Assume $T_1<T_2$. The origin is an equilibrium point of the system so $x(T_2;u_1)^+=0$ and $J_{T_2}({u_1}_{|[0;T_2]})=J_{T_1}({u_1}_{|[0;T_1]})=J_{T_2}({u_2}_{|[0;T_2]})$.  Applying Lemma~2 in the hyperbolic case, the optimal control for the finite horizon problem in time $T_2$ is unique. Then ${u_2}_{|[0;T_2]}={u_1}_{|[0;T_2]}$ and $u_1=u_2$ in $\Ui$. $\square$

\begin{corollary} If $A$ is hyperbolic, then $\mu_\infty:\A_\infty\to\R$ is a continuous function.
\end{corollary}
\noindent Indeed, according to Theorem~3, optimal control has compact support. This support depends on $x_0$ but we see from the proof of the theorem that it remains uniformly bounded if $x_0$ runs over a compact subset of $\A_\infty$. Continuity of
$\mu_\infty$ now follows from Theorem~2.

\begin{theorem}
     Assume that $E^+=\R^n$. A control $u\in\Ui$ such that $x(T;u)=0$ and $x(t;u)\neq0$ for any $t<T$ is optimal for the infinite horizon problem if and only if $x(.;u)$ can be complemented by $p(\cdot)$ in such a way that the Pontryagin maximum principle is satisfied, and, moreover, $|p(T)B|=1$ and $|p(T)e^{-\tau A}B|\leq 1, \forall \tau \geq 0$.
\end{theorem}

\noindent{\bf Proof.} It follows from Theorem~3 that optimal control is also optimal for the problem with free finite time. In this case, the transversality condition states that $h_{u(t)}(p(t),x(t;u)) = p(t)Ax(t;u) + p(t)Bu(t) - |u(t)| = 0, \forall t \geq 0$. There exists $\varepsilon > 0$ such that $u(t) = \pm 1, T-\varepsilon \leq t \leq T$, otherwise $x(T-\varepsilon;u)=0$. Then $p(t)Bu(t) - |u(t)| = |p(t)B| - 1, T - \varepsilon \leq t \leq T$ (maximality condition), and because $x(T;u)=0$, we obtain $|p(T)B| = 1$.

Assume that $u$ is optimal for the free time problem and that $x(T;u)=0$, we may assume that $u(t) = 0, \forall t \geq T$. Then the control $u(t), 0 \leq t \leq T+s$ is optimal for the free time problem for any $s>0$. Hence $ x(t;u),\ 0 \leq t \leq T+s $, can be complemented by $p_s(\cdot)$ in such a way that the Pontryagin maximum principle is satisfied. Then $|p_s(T)B|=1$ and $|p_s(T+ \tau)B| = |p_s(T)e^{- \tau A}B| \leq 1,\: \forall\,\tau\in[0,s]$.

\begin{lemma}
    Let $s>0$; then the set $\{ \xi \in \mathbf{R}^{n*}: |\xi e^{-\tau A}B| \leq 1, 0\leq \tau \leq s \}$ is bounded.
\end{lemma}
\noindent{\bf Proof of the lemma.} By contradiction assume that there exists a sequence $(\xi_k)$ in the previous set such that $|\xi_k| \to +\infty$ as $k \to \infty$. We may assume that $\frac{\xi_k}{|\xi_k|} \to \eta \in \mathbf{R}^{n*}$, $|\eta|=1$. Then passing to the limit for $0 \leq \tau \leq s$, we have $\eta e^{-\tau A}B = 0$, which contradicts the Kalman condition. $\square$ \\

Applying this lemma, the family $\{p_s(T), s>0\}$ is bounded so it has limiting points as $s \to \infty$. Any limiting point satisfies the conditions of the theorem. Hence these are necessary conditions for optimality.

Moreover, the Pontryagin maximum principle is a sufficient optimality condition for the fixed time problem. It follows that our extremal is optimal for every time $T+s, s\geq 0$, i.e. for any arbitrarily fixed time. Hence this extremal is optimal for the free time problem and conditions of the theorem are sufficient for optimality. $\square$

\section{Two-dimensional Systems}

\subsection{Classical examples}

In this section we study the stabilization to the origin in infinite time horizon of classical two-dimensional systems with scalar control $|u| \leq 1$: the case of a free particle and the one of a harmonic oscillator.

\subsubsection{Free Particle}

Considering the degenerate case of a free particle, we study the stabilization of the 2-dimensional system to $0_{\mathbf{R}^2}$ in infinite time horizon:

\begin{equation*}
    \begin{pmatrix}
        \dot{x}^1 \\
        \dot{x}^2
    \end{pmatrix} = \begin{pmatrix}
        0&1 \\
        0&0
    \end{pmatrix} \begin{pmatrix}
        x^1 \\
        x^2
    \end{pmatrix} + u \begin{pmatrix}
        0 \\
        1
    \end{pmatrix} \qquad |u| \leq 1.
\end{equation*}

\begin{theorem}
    For every $x_0 \in \mathbf{R}^2$, $\mu_{\infty}(x_0)=\inf\limits_{T>0}{\mu_T(x_0)} = |x^2_{0}|$. \\
    If $x_0 \in \{x \: | \: x^1 \leq -\frac{1}{2}(x^2)^2,x^2\geq 0 \} $ or $x \in \{x \: | \: x^1 \geq \frac{1}{2}(x^2)^2,x^2\leq 0 \}$, then there exists $0<T<\infty$ such that $\mu_T(x_0) = |x^2_{0}|$. Otherwise $\mu_{T}(x_0) > |x^2_{0}|$ for all $0<T<\infty$.
\end{theorem}

\begin{center}
\includegraphics[scale=0.35]{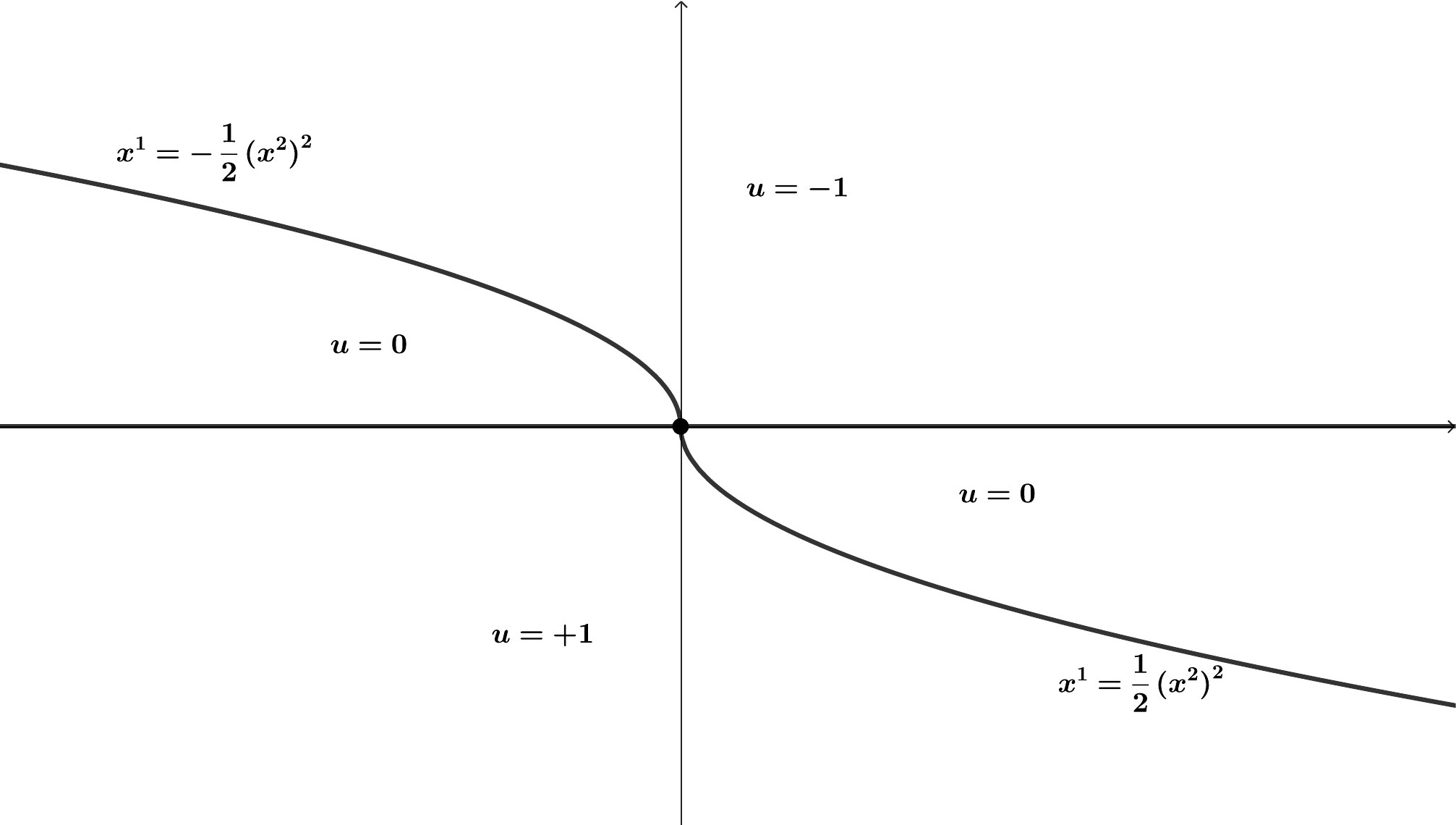}
\captionof{figure}{Regions of optimal control at the limit in infinite time horizon}
\end{center}

\noindent{\bf Proof.}
Noting that $\dot{x}^2 = u$, for every $x_{0} \in \mathbf{R}^2$, for every $T>0$:
\begin{equation*}
    \mu_T(x_0) = \int_0^T|u(t)|dt \geq |\int_0^Tu(t)dt| = |\int_0^T\dot{x}^2(t)dt| = |x^2(0) - x^2(T)| = |x^2_{0}|.
\end{equation*}
Then $\mu_{\infty}(x_0)=\inf\limits_{T>0}{\mu_T(x_0)} \geq |x^2_{0}|$. \\
We study the form of the trajectory in function of the control. \\
If $u \equiv 1$ or $u \equiv -1$:
\begin{equation*}
    \begin{cases}
        \dot{x}^1 = x^2 \\
        \dot{x}^2 = \pm 1
    \end{cases}, \qquad x_1 = \pm \frac{1}{2}((x^2)^2 - (x^2_{0})^2)+x^1_{0}.
\end{equation*}
If $u\equiv0$:
\begin{equation*}
    \begin{cases}
        \dot{x}^1 = x^2 \\
        \dot{x}^2 = 0
    \end{cases}, \qquad \begin{cases}
        x^1(t) = x^2_{0}t + x^1_{0} \\
        x^2(t) = x^2_{0}.
    \end{cases}
\end{equation*}
\\
\underline{Region where $\mu_{\infty}(x_0)$ in reached in finite time:} \\
If $x^1_{0} = \frac{1}{2}(x^2_{0})^2$ and $x^2_{0}\leq0$ we can take a control $u \equiv +1$ until $x$ reaches the origin (and then $u \equiv 0)$. With this control, we obtain $\mu_T(x_0) = |x^2_{0}|$ for every $T \geq |x^2_{0}|$. By symmetry, if $x^1_{0} = -\frac{1}{2}(x^2_{0})^2$ and $x^2_{0}\geq0$, with a control $u\equiv-1$ until the origin (then $u \equiv 0$), we obtain the same result.
\\ \\
If now $x^1_{0} > \frac{1}{2}(x^2_{0})^2$ and $x^2_{0}\leq0$, we can take a control $u \equiv 0$ until $x^1 = \frac{1}{2}(x^2)^2$. The arc of the trajectory is a horizontal line of equation $x^2 = x^2_{0}$. Then we take $u \equiv +1$ and the trajectory is the same as the previous one. The cost does not change by adding a time interval where $u \equiv 0$. By symmetry, we obtain the same result when $x^1_{0} < -\frac{1}{2}(x^2_{0})^2$ and $x^2_{0}\geq0$.

\begin{center}
\includegraphics[scale=0.35]{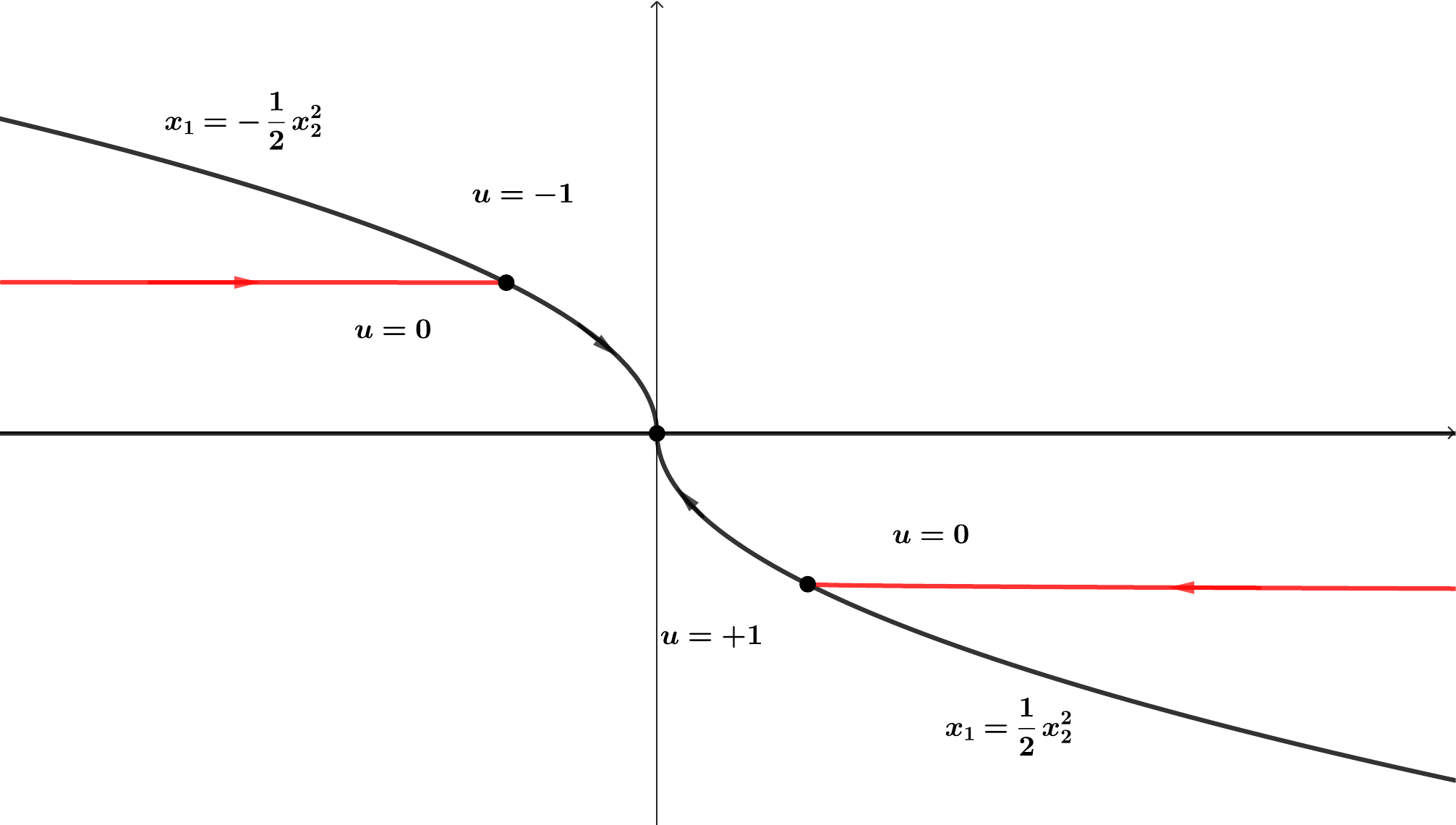}
\captionof{figure}{Optimal trajectories in finite time}
\end{center}
\underline{Region where $\mu_{\infty}(x_0)$ is not reached in finite time:} \\
We construct a sequence $T_k$ such that $\mu_{T_k}(x_0) \rightarrow |x^2_{0}|$. Applying the PMP, we obtain some extremals, that are all optimal thank to the linearity of the system and the convexity of the set of constraints. Solving the adjoint equation $\dot{p}=-pA$, we obtain the switching function:
\begin{equation*}
    f(t) = p(t)B = -p^1_{0}t + p^2_{0}.
\end{equation*}
Variations of the co-vector allows the switches $-1 \rightarrow 0 \rightarrow +1$  and $+1 \rightarrow 0 \rightarrow -1$ with every lengths for the different time intervals where $u$ is constant.
\\ \\
If $x^1_{0} > -\frac{1}{2}(x^2_{0})^2$ and $x^2_{0}\geq0$, we consider a control with switches $-1 \rightarrow 0 \rightarrow +1$. First arc of the trajectory is parabolic with equation $ x^1 = -\frac{1}{2}((x^2)^2 - (x^2_{0})^2)+x^1_{0}$. We choose the first time of switch when $x^2(t) = -\frac{1}{k}$. Second arc is then a horizontal line of equation $x^2=-\frac{1}{k}$. Third arc is parabolic with equation $x^1=+\frac{1}{2}(x^2)^2$ and reaches the origin at a time denoted $T_k$. This control satisfies the equation of the PMP, thus is optimal. \\
Because $x^2_{0}\geq0$ we have:
\begin{equation*}
    \mu_{T_k}(x_0)=|x^2_{0}-(-\frac{1}{k})| + |\frac{1}{k} - 0 | = x^2_{0} + \frac{2}{k} \underset{k\rightarrow + \infty}{\longrightarrow} x^2_{0}.
\end{equation*}
So $\mu_{\infty}(x_0) = \inf\limits_{T>0}\mu_T(x_0) = |x^2_{0}|$.
If $x^1_{0} < \frac{1}{2}(x^2_{0})^2$ and $x^2_{0}\leq0$ by symmetry we can construct such a sequence with a control $+1 \rightarrow 0 \rightarrow -1$ and we obtain also $\mu_{\infty}(x_0) = |x^2_{0}|$.
\begin{center}
\includegraphics[scale=0.45]{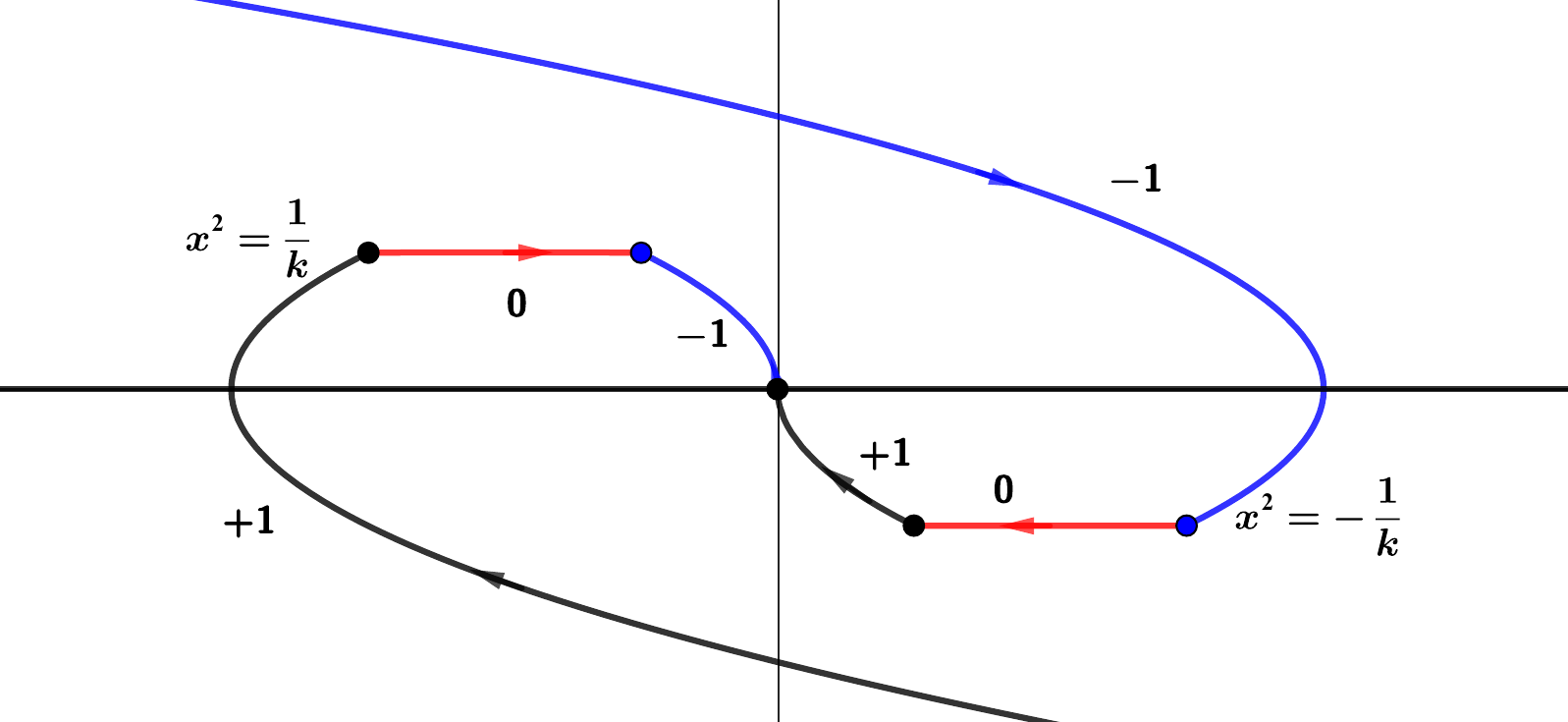}
\captionof{figure}{Construction of an optimal trajectory}
\end{center}
The trajectory limit when $k$ goes to infinity is not admissible, because if $x^2 = 0$, $x^1 \neq 0$ and $u \equiv 0$, then the point is an equilibrium point of the system and cannot reach the origin. We also see that as soon as the optimal control changes sign, then $\mu_T(x_0) > |x^2_{0}|$. $\square$

\subsubsection{Harmonic Oscillator}

Considering the example of a 2-dimensional controlled harmonic oscillator, we study the stabilization to the origin in infinite time horizon:
\begin{equation*}
    \begin{pmatrix}
        \dot{x}^1 \\
        \dot{x}^2
    \end{pmatrix} = \begin{pmatrix}
        0&1 \\
        -1&0
    \end{pmatrix} \begin{pmatrix}
        x^1 \\
        x^2
    \end{pmatrix} + u \begin{pmatrix}
        0 \\
        1
    \end{pmatrix}.
\end{equation*}

\begin{prop}
    In a finite time $T>0$, the optimal control $u$ is unique and has a periodic structure after some first switching time $\alpha_0 \in [0;\pi)$. Let $\varepsilon \in \{-1;+1\}$. There exists a first switching time $\alpha_0\in [0;\pi)$ and a period $\delta \in [0; \pi)$ such that, for $t \geq \alpha_0$, $u$ takes the values $-1$, $0$ and $+1$ in a periodic way :
    \begin{equation*}
        -\varepsilon \rightarrow 0 \rightarrow +\varepsilon \rightarrow 0 \rightarrow -\varepsilon \rightarrow ... \qquad or \qquad 0 \rightarrow +\varepsilon \rightarrow 0 \rightarrow -\varepsilon \rightarrow 0 \rightarrow ...
    \end{equation*}
The time where $u$ is $+1$ or $-1$ is denoted $\delta$ and then the one where $u$ is $0$ is equal to $\pi - \delta$. So the optimal control $u$ takes one of the following forms until the final time $T$:
\begin{equation*}
    u(t)=\begin{cases}
        0 \: \: \: \: \text{if} \: \: t \in [0,\alpha_0)\cup_{k\geq 0}[\alpha_0 + k\pi + \delta, \alpha_0 +(k+1)\pi) \\
        \pm \varepsilon \: \text{if} \: \: t \in \cup_{k \geq 0} [\alpha_0 + 2k\pi, \alpha_0 + 2k\pi + \delta) \\
        \mp \varepsilon \: \text{if} \: \: t \in \cup_{k \geq 0} [\alpha_0 + (2k+1)\pi,\alpha_0 + (2k+1)\pi + \delta)
    \end{cases} \qquad \alpha_0 \leq \pi - \delta,
\end{equation*}
or
\begin{equation*}
    u(t)=\begin{cases}
        \pm \varepsilon \: \: \: \: \text{if} \: \: t \in [0,\alpha_0) \cup_{k \geq 0} [\alpha_0 + 2k\pi - \delta, \alpha_0 + 2k\pi) \\
        \mp \varepsilon \: \: \: \: \text{if} \: \: t \in \cup_{k \geq 0}[\alpha_0 + (2k+1)\pi -\delta, \alpha_0 + (2k+1)\pi) \\
        0 \: \: \: \: \text{if} \: \: t \in \cup_{k \geq 0} [\alpha_0+ k\pi, \alpha_0 + (k+1)\pi - \delta)
    \end{cases} \qquad \alpha_0 \leq \delta.
\end{equation*}

\end{prop}

\noindent{\bf Proof.} Applying Lemma~1, $B$ is a column and $\det(A) \neq 0$, the PMP gives the description of the optimal control, which is unique. Solving the adjoint equation $\dot{p}=-pA$, we obtain the switching function:
\begin{equation*}
    f(t) = p(t)B = a \sin(t+b),
\end{equation*}
where $a$ and $b$ are real parameters related to $p^1_{0},p^2_{0}$. Then if $|a|>1$, the optimal control is not always $0$ and takes the structure described in the proposition. $\square$

\begin{prop}
The optimal trajectory alternates periodically between different arc of circles, traveled clockwise, as the optimal control changes sign: \\
Circles centered at $(\pm 1,0)$ when $u = \pm 1$, \\
Circle centered at $(0,0)$ when $u=0$.
\end{prop}

\noindent{\bf Proof.} If $u \equiv 0$:
\begin{equation*}
    \begin{cases}
          \dot{x}^1(t) = x^2(t) \\
          \dot{x}^2(t) = - x^1(t)
    \end{cases}, \qquad \begin{cases}
        x^1(t) = a_0\cos(t+b_0) \\
        x^2(t) = a_0\sin(t+b_0).
    \end{cases}
\end{equation*}
Then:
\begin{equation*}
    \frac{d}{dt} \{ (x^1)^2(t) + (x^2)^2(t) \} = 0,
\end{equation*}
so the trajectory is a clockwise circle centered at the origin. \\
If $u \equiv 1$:
\begin{equation*}
    \begin{cases}
          \dot{x}^1(t) = x^2(t) \\
          \dot{x}^2(t) = - x^1(t) + 1
    \end{cases}, \qquad \begin{cases}
        x^1(t) = a_1\cos(t+b_1) + 1  \\
          x^2(t) = a_1\sin(t+b_1).
    \end{cases}
\end{equation*}
Then:
\begin{equation*}
    \frac{d}{dt} \{ (x^1(t)-1)^2 + (x^2)^2(t) \} = 0,
\end{equation*}
so the trajectory is a clockwise circle centered at $(1,0)$. \\
If $u \equiv -1$, by symmetry the trajectory is a clockwise circle centered at $(-1,0)$. $\square$

\begin{theorem}
    Let $x_0 \in \mathbf{R}^2$. For every $k \in \mathbf{N}$ such that $|k| \geq |x_0|$, there exists a finite time $T(k)$, $(k-1)\pi \leq T(k) \leq k\pi$, such that the optimal trajectory which returns to the origin has its switches on the circles centered at $(0,n)$ for $n \in \mathbf{Z}$, $n \leq k$. The optimal cost of this trajectory is $c(k) = k \arccos(1 - \frac{|x_0|^2}{2k^2})$. The limit of the optimal cost is $\mu_{\infty}(x_0)=|x_0|$ and the point-wise limit of the optimal trajectory is the clockwise circle centered at the origin $x_{\infty}(t) = \begin{pmatrix}
        \cos(t) & \sin(t) \\
        - \sin(t) & \cos(t)
    \end{pmatrix}x_0$.
\end{theorem}

\noindent{\bf Proof.} The proof is presented in Appendix~A. 

\begin{center}
\includegraphics[scale=0.35]{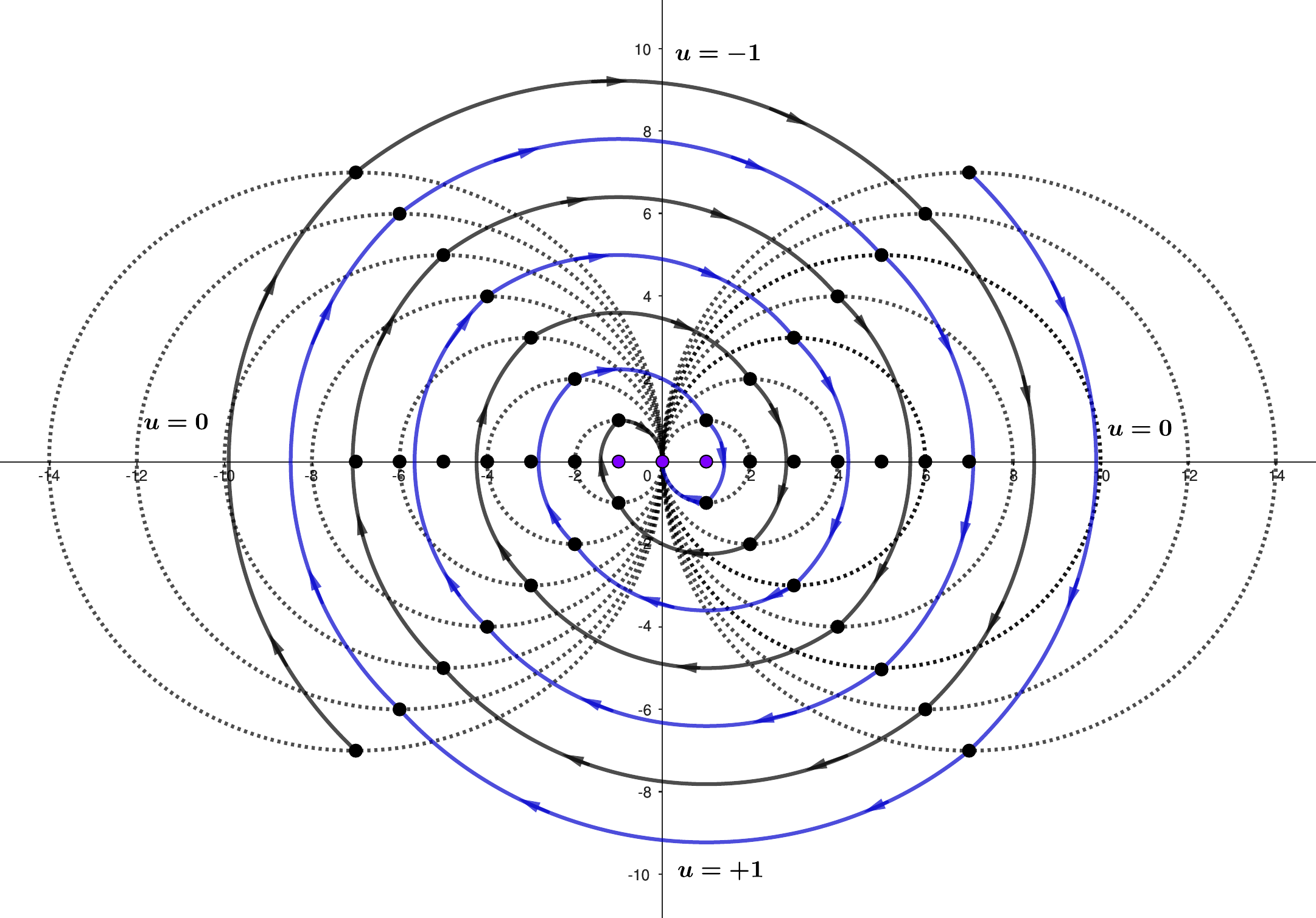}
\captionof{figure}{Switching curves and optimal trajectories for the harmonic oscillator}
\end{center}

\subsection{Hyperbolic systems}

In this section, we study the stabilization to the origin in infinite time of all real hyperbolic two-dimensional systems $\dot{y} = Ay + Bu$ with a scalar control $|u| \leq 1$. The matrix A is similar to one of the following where $(\lambda_i,\lambda,\alpha,\beta) \in \mathbf{R}^4$:
\begin{equation*}
    \begin{pmatrix}
        \lambda_1 & 0 \\
        0 & \lambda_2
    \end{pmatrix}, \qquad
    \begin{pmatrix}
        \lambda & 1 \\
        0 & \lambda
    \end{pmatrix}, \qquad
    \begin{pmatrix}
        \alpha & \beta \\
        - \beta & \alpha
    \end{pmatrix}.
\end{equation*}
As explained in Theorem~3, $\mu_{\infty}(x_0)$ only depends on the part $x_0^+$ corresponding to the eigenvalues of $A$ with positive real parts. Thus we only consider cases where $\mathbf{R}^2 = E^+$, i.e. where $x_0 = x_0^+$ for all $x_0 \in \mathbf{R}^2$ : $\lambda_i,\lambda,\alpha > 0$.
In order to treat all cases $B = \begin{pmatrix}
    b_1 \\
    b_2
\end{pmatrix}$ with $(A,B)$ respecting the Kalman condition, we apply a change of basis $x = Py$ such that $P$ and $A$ commute. Thus the system in these new coordinates can be written $\dot{x} = Ax + PBu$.
\subsubsection{First Hyperbolic case : $A = \begin{pmatrix}
    \lambda_1 & 0 \\
    0 & \lambda_2
\end{pmatrix}$}

The Kalman condition is respected if and only if $b_1, b_2 \neq 0$. Applying the change of basis $P = \begin{pmatrix}
b_1^{-1} & 0 \\
0 & b_2^{-1}
\end{pmatrix}$, the study reduces to the stabilization of the system:
\begin{equation*}
    \begin{pmatrix}
        \dot{x}^1 \\
        \dot{x}^2
    \end{pmatrix} = \begin{pmatrix}
        \lambda_1 & 0 \\
        0 & \lambda_2
    \end{pmatrix} \begin{pmatrix}
        x^1 \\
        x^2
    \end{pmatrix} + u \begin{pmatrix}
        1 \\
        1
    \end{pmatrix} \qquad \lambda_1, \lambda_2>0 \qquad |u| \leq 1.
\end{equation*}

\begin{prop}
    Let $x_0 \in \mathbf{R}^2$. If $0_{\mathbf{R}^2}$ is reachable from $x_0$, then the optimal control is unique and has isolated switches with one of the following structure:
    \begin{equation*}
        \varepsilon \rightarrow 0, \qquad
        0 \rightarrow \varepsilon \rightarrow 0, \qquad
        - \varepsilon \rightarrow 0 \rightarrow \varepsilon \rightarrow 0 \qquad \varepsilon \in \{-1;+1\}
    \end{equation*}
    The optimal control $u$ has thus one of the following forms for $\varepsilon \in \{-1;+1\}$:
    \begin{equation*}
        u(t)=\begin{cases} \varepsilon \: \: \: \: \text{for} \: \: t \in [0,t_1) \\
        0 \: \: \: \: \text{for} \: \: t \in [t_1,+\infty)
        \end{cases} \qquad t_1 \geq 0,
    \end{equation*}
    or 
    \begin{equation*}
        u(t)=\begin{cases}
            0 \: \: \: \: \text{for} \: \: t \in [0,t_1)\cup[t_2,+\infty) \\
            \varepsilon \: \: \: \: \text{for} \: \: t \in [t_1,t_2)
        \end{cases} \qquad 0 \leq t_1 \leq t_2,
    \end{equation*}
    or
    \begin{equation*}
        u(t)=\begin{cases}
            - \varepsilon \: \text{for} \: \: t \in [0,t_1) \\
            0 \: \: \: \: \text{for} \: \: t \in [t_1,t_2)\cup[t_3,+\infty) \\
            \varepsilon \: \: \: \: \text{for} \: \: t \in[t_2,t_3) 
        \end{cases} \qquad 0 \leq t_1 \leq t_2 \leq t_3.
    \end{equation*}
\end{prop}

\noindent{\bf Proof.} The uniqueness and the structure of the optimal control $u$ are given applying Lemma~2, noting that A is hyperbolic. Solving the adjoint equation $\dot{p} = -p A$, we obtain the switching function:
\begin{equation*}
    f(t) = p(t)\cdot B = p^1_{0}e^{-\lambda_1 t} + p^2_{0}e^{-\lambda_2 t}.
\end{equation*}
If $f(t)$ is less than $-1$, then $u \equiv -1$. If it is more than $+1$, then $u \equiv +1$. Else $u \equiv 0$. The function goes to 0 at infinity and has at most two different monotone branches. Thus we obtain one of the structures described in the proposition. $\square$

\begin{theorem}
    The region of initial conditions from which there exists a
    trajectory that reaches $0_{\mathbf{R}^2}$ is open and delimited by arcs of the curves $\mathcal{C}_+$ and $\mathcal{C}_-$ between the points $(\frac{1}{\lambda_1},\frac{1}{\lambda_2})$ and $(-\frac{1}{\lambda_1},-\frac{1}{\lambda_2})$:
    \begin{equation*}
        \mathcal{C}_+ \: : \: x^2 = \frac{1}{\lambda_2}(2(\frac{\lambda_1 x^1 + 1}{2})^{\frac{\lambda_2}{\lambda_1}}-1), \qquad  \mathcal{C}_- \: : \: x^2 = \frac{1}{\lambda_2}(1-2(\frac{1-\lambda_1 x^1}{2})^{\frac{\lambda_2}{\lambda_1}}).
    \end{equation*}
    The switching curve $0\rightarrow+1$ is given by the arc of the curve $\mathcal{C}_{0,+1}$ between the points $(0,0)$ and $(-\frac{1}{\lambda_1},-\frac{1}{\lambda_2})$ (respectively $\mathcal{C}_{0,-1}$, $(0,0)$, $(\frac{1}{\lambda_1},\frac{1}{\lambda_2})$ for the switch $0 \rightarrow -1$):
    \begin{equation*}
        \mathcal{C}_{0,+1} \: : \: x^2 = \frac{1}{\lambda_2}((\lambda_1 x^1 + 1)^{\frac{\lambda_2}{\lambda_1}}-1), \qquad \mathcal{C}_{0,-1} \: : \: x^2 = \frac{1}{\lambda_2}(-(1-\lambda_1 x^1)^{\frac{\lambda_2}{\lambda_1}}+1).
    \end{equation*}
    The switching curve $+1 \rightarrow 0$ is given by the arc of the curve $\mathcal{C}_{+1,0}$ between points $(0,0)$ and $(\frac{1}{\lambda_1},\frac{1}{\lambda_2})$ (respectively $\mathcal{C}_{-1,0}$, $(0,0)$ and $(-\frac{1}{\lambda_1},-\frac{1}{\lambda_2})$ for the switch $-1 \rightarrow +1$):
    \begin{equation*}
        \mathcal{C}_{+1,0} \: : \: x^2 = \frac{1}{\lambda_2} (\lambda_1 x^1)^{\frac{\lambda_2}{\lambda_1}}, \qquad \mathcal{C}_{-1,0} \: : \: x^2 = -\frac{1}{\lambda_2}(-\lambda_1 x^1)^{\frac{\lambda_2}{\lambda_1}}.
    \end{equation*}
\end{theorem}

\begin{center}
\includegraphics[scale=0.45]{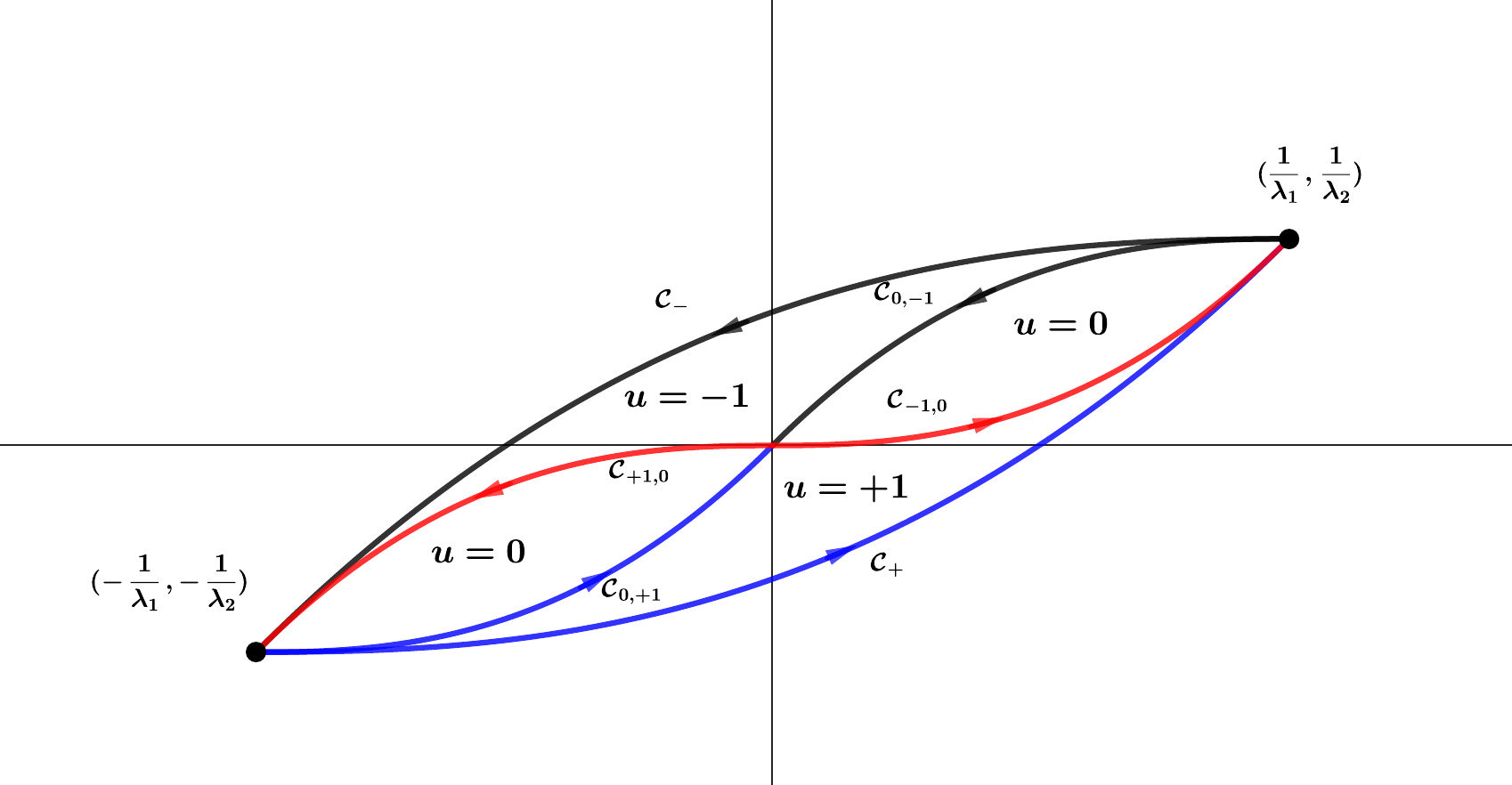}
\captionof{figure}{Attainable set and switching curves in the case $\lambda_1 < \lambda_2$}
\end{center}

\noindent{\bf Proof.} First we look at the form of the trajectories when the control is constant. If $u\equiv+1$:
\begin{equation*}
    \begin{cases}
        \dot{x}^1 = \lambda_1 x^1 + 1 \\
        \dot{x}^2 = \lambda_2 x^2 + 1
    \end{cases},
\end{equation*}
the only equilibrium point is $(-\frac{1}{\lambda_1},-\frac{1}{\lambda_2})$. Otherwise the solution $|x(t)|$ goes to $+\infty$ in infinite time. Solving the differential system from an initial condition $(x^1_{0},x^2_{0})$ we obtain a trajectory of the following form, traveled towards the infinity:
\begin{equation*}
    \frac{\lambda_2 x^2 +1}{\lambda_2 x^2_{0} + 1} = (\frac{\lambda_1 x^1 + 1}{\lambda_1 x^1_{0} + 1})^{\frac{\lambda_2}{\lambda_1}}.
\end{equation*}
If $(0,0)$ belongs to the trajectory then $(x^1_{0},x^2_{0})$ belongs to the arc of the curve $\mathcal{C}_{0,+1}$ between $(-\frac{1}{\lambda_1},-\frac{1}{\lambda_2})$ and $(0,0)$.
\\ \\
Then we look at the trajectories with one switch of type $0 \rightarrow +1$. If $u\equiv0$ on a time interval:
\begin{equation*}
    \begin{cases}
        \dot{x}^1 = \lambda_1 x^1 \\
        \dot{x}^2 = \lambda_2 x^2
    \end{cases},
\end{equation*}
the only equilibrium point is $(0,0)$. Otherwise the trajectory has the following form, traveled towards the infinity:
\begin{equation*}
    \frac{x^2}{x^2_{0}}=(\frac{x^1}{x^1_{0}})^{\frac{\lambda_2}{\lambda_1}}.
\end{equation*}
In order to arrive at $0_{\mathbf{R}^2}$ with one switch $0 \rightarrow +1$, the trajectory must cross the arc of $\mathcal{C}_{0,+1}$ between the points $(-\frac{1}{\lambda_1},-\frac{1}{\lambda_2})$ and $(0,0)$. Because trajectories corresponding to $u \equiv 0$ do not cross the curve except at the origin, the region of initial conditions for which such a switch is possible is delimited by the trajectory with control $0$ which passes through the point $(-\frac{1}{\lambda_1},-\frac{1}{\lambda_2})$. We obtain the arc of the curve $\mathcal{C}_{-1,0}$ between the points $(0,0)$ and $(-\frac{1}{\lambda_1},-\frac{1}{\lambda_2})$.
\\ \\
By symmetry we obtain the switching curves for $0 \rightarrow -1$ between the points $(0,0)$ and $(\frac{1}{\lambda_1},\frac{1}{\lambda_2})$. Finally in order to obtain a switch of type $-1 \rightarrow 0 \rightarrow +1$, the trajectory with $u\equiv-1$ at the beginning must cross the arc of $\mathcal{C}_{-1,0}$ between the points $(0,0)$ and $(-\frac{1}{\lambda_1},-\frac{1}{\lambda_2})$. Because the trajectories do not cross except in $(\frac{1}{\lambda_1},\frac{1}{\lambda_2})$ when $u\equiv-1$, the region of initial conditions for which such a switch is possible is delimited by the one which passes through the point $(-\frac{1}{\lambda_1},-\frac{1}{\lambda_2})$. Thus we obtain the arc of the curve $\mathcal{C}_{-}$ between the points $(\frac{1}{\lambda_1},\frac{1}{\lambda_2})$ and $(-\frac{1}{\lambda_1},-\frac{1}{\lambda_2})$. $\square$

\begin{center}
\includegraphics[scale=0.45]{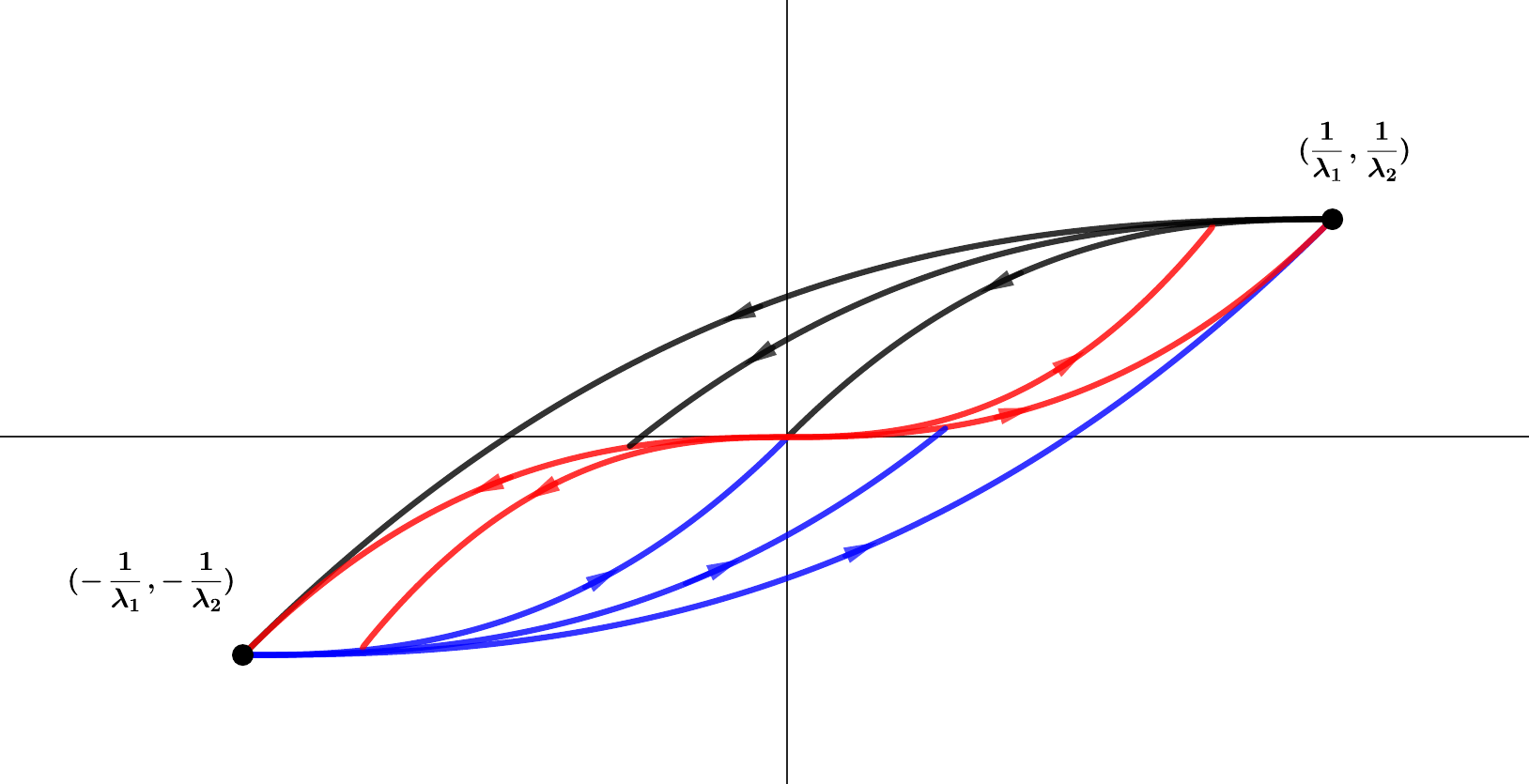}
\captionof{figure}{Optimal trajectories}
\end{center}

\subsubsection{Second Hyperbolic case : $A = \begin{pmatrix}
    \lambda & 1 \\
    0 & \lambda
\end{pmatrix}$}

The Kalman condition is respected if and only if $b_2 \neq 0$. Applying the change of basis $P = \begin{pmatrix}
b_2^{-1} & 0 \\
0 & b_2^{-1}
\end{pmatrix}$ and denoting $b = \frac{b_1}{b_2}$, the study reduces to the stabilization of the system:
\begin{equation*}
    \begin{pmatrix}
        \dot{x}^1 \\
        \dot{x}^2
    \end{pmatrix} = \begin{pmatrix}
        \lambda & 1 \\
        0 & \lambda
    \end{pmatrix} \begin{pmatrix}
        x^1 \\
        x^2
    \end{pmatrix} + u\begin{pmatrix}
        b\\
        1
    \end{pmatrix} \qquad \lambda > 0 \qquad |u| \leq 1.
\end{equation*}

\begin{prop}
    Let $x_0 \in \mathbf{R}^2$. If $0_{\mathbf{R}^2}$ is reachable from $x_0$, then the optimal control is unique and has isolated switches with one of the following forms for $\varepsilon \in \{-1;+1\}$:
    \begin{equation*}
        u(t)=\begin{cases} \varepsilon \: \: \: \: \text{for} \: \: t \in [0,t_1) \\
        0 \: \: \: \: \text{for} \: \: t \in [t_1,+\infty)
        \end{cases} \qquad t_1 \geq 0,
    \end{equation*}
    or 
    \begin{equation*}
        u(t)=\begin{cases}
            0 \: \: \: \: \text{for} \: \: t \in [0,t_1)\cup[t_2,+\infty) \\
            \varepsilon \: \: \: \: \text{for} \: \: t \in [t_1,t_2)
        \end{cases} \qquad 0 \leq t_1 \leq t_2,
    \end{equation*}
    or
    \begin{equation*}
        u(t)=\begin{cases}
            - \varepsilon \: \text{for} \: \: t \in [0,t_1) \\
            0 \: \: \: \: \text{for} \: \: t \in [t_1,t_2)\cup[t_3,+\infty) \\
            \varepsilon \: \: \: \: \text{for} \: \: t \in[t_2,t_3) 
        \end{cases} \qquad 0 \leq t_1 \leq t_2 \leq t_3.
    \end{equation*}
\end{prop}

\noindent{\bf Proof.} The proof is the same as in the previous study of the first hyperbolic case, we apply the Lemma~2 and we note that the switching function, which is given by the following formula, has at most two monotone branches and goes to 0 at infinity:
\begin{equation*}
    f(t) = p(t) \cdot B = (p^1_{0} + p^2_{0} - p^2_{0}t)e^{-\lambda t}. \qquad \square
\end{equation*}

\begin{theorem}
    The region of initial conditions from which there exists a
    trajectory that reaches $0_{\mathbf{R}^2}$ is open and delimited by arcs of the curves $\mathcal{C}_+$ and $\mathcal{C}_-$ between the points $(\frac{1}{\lambda}(b-\frac{1}{\lambda}),\frac{1}{\lambda})$ and $(-\frac{1}{\lambda}(b-\frac{1}{\lambda}),-\frac{1}{\lambda})$:
    \begin{align*}
        \mathcal{C}_+ \: &: \: x^1 = -\frac{1}{\lambda}(b-\frac{1}{\lambda}) + \frac{1}{\lambda^2}(\lambda x^2 + 1)(\lambda b - 1 + \ln(\frac{\lambda x^2 + 1}{2})), \\
        \mathcal{C}_- \: &: \: x^1 = \frac{1}{\lambda}(b-\frac{1}{\lambda}) + \frac{1}{\lambda^2}(\lambda x^2 - 1)(\lambda b - 1 + \ln(\frac{1-\lambda x^2}{2})).    \end{align*}
    The switching curve $0\rightarrow+1$ is given by the arc of the curve $\mathcal{C}_{0,+1}$ between the points $(0,0)$ and $(-\frac{1}{\lambda}(b-\frac{1}{\lambda}),-\frac{1}{\lambda})$ (respectively $\mathcal{C}_{0,-1}$, $(0,0)$, $(\frac{1}{\lambda}(b-\frac{1}{\lambda}),\frac{1}{\lambda})$ for the switch $0 \rightarrow -1$):
    \begin{align*}
        \mathcal{C}_{0,+1} \: &: \: x^1 = -\frac{1}{\lambda}(b-\frac{1}{\lambda}) + \frac{1}{\lambda^2}(\lambda x^2 + 1)(\lambda b - 1 + \ln(\lambda x^2 + 1)), \\
        \mathcal{C}_{0,-1} \: &:  \: x^1 = \frac{1}{\lambda}(b-\frac{1}{\lambda}) + \frac{1}{\lambda^2}(\lambda x^2 - 1)(\lambda b - 1 + \ln(1-\lambda x^2)).
    \end{align*}
    The switching curve $+1 \rightarrow 0$ is given by the arc of the curve $\mathcal{C}_{+1,0}$ between the points $(0,0)$ and $(\frac{1}{\lambda}(b-\frac{1}{\lambda}),\frac{1}{\lambda})$ (and respectively $\mathcal{C}_{-1,0}$, $(0,0)$ and $(-\frac{1}{\lambda}(b-\frac{1}{\lambda}),-\frac{1}{\lambda})$ for the switch $-1 \rightarrow +1$):
    \begin{equation*}
        \mathcal{C}_{+1,0} \: : \: x^1 = x^2(1-\frac{1}{\lambda} + \frac{1}{\lambda}\ln(-\lambda x^2)), \qquad \mathcal{C}_{-1,0} \: :  \: x^1 = x^2(1-\frac{1}{\lambda} + \frac{1}{\lambda}\ln(\lambda x^2)).
    \end{equation*}
\end{theorem}

\begin{center}
\includegraphics[scale=0.35]{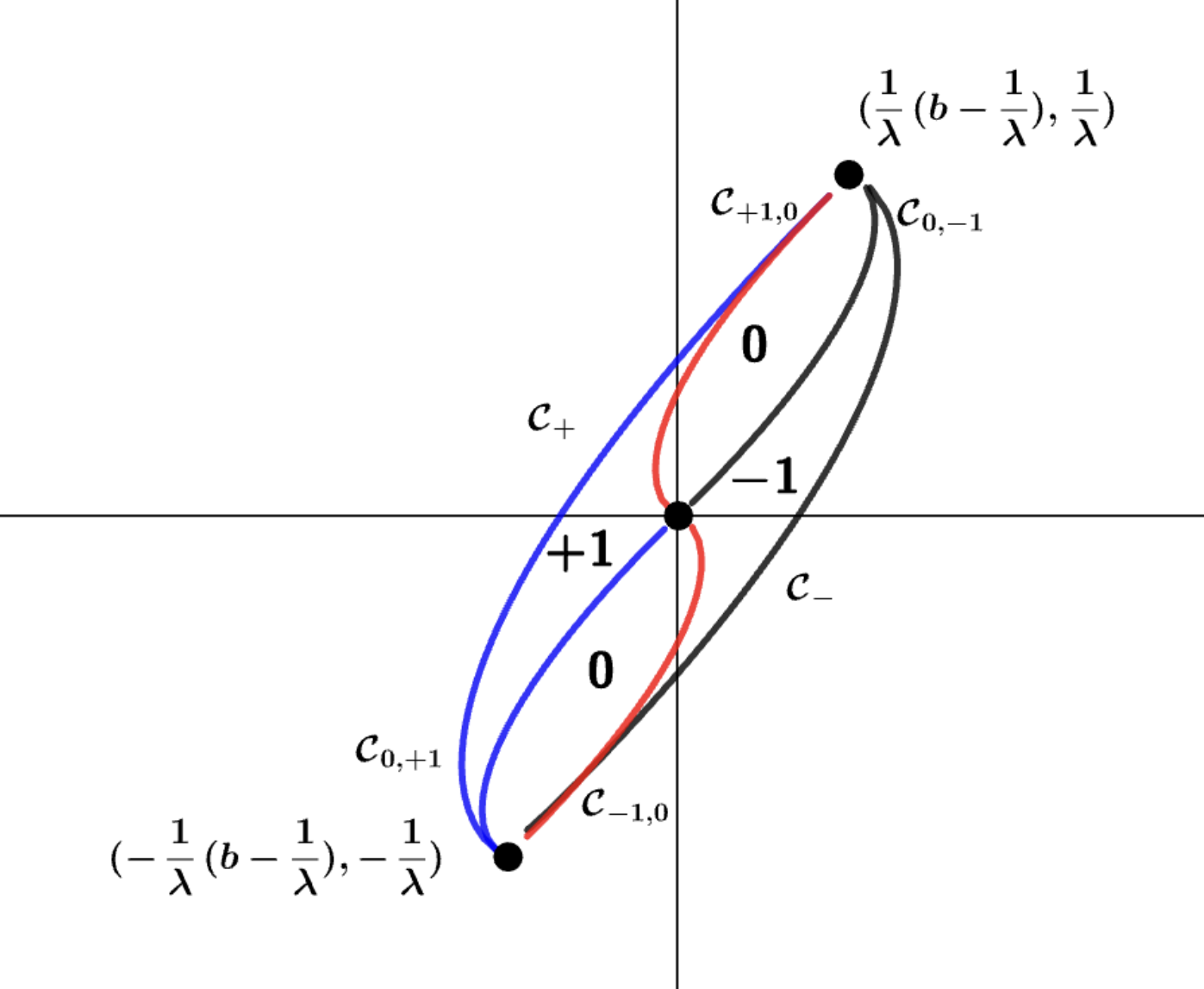}
\captionof{figure}{Attainable set and switching curves}
\end{center}

\noindent{\bf Proof.} The scheme of the study is exactly the same as in the first hyperbolic case. 

\subsubsection{Third Hyperbolic case : $A = \begin{pmatrix}
    \alpha & \beta \\
    - \beta & \alpha
\end{pmatrix}$}

The Kalman condition is respected if and only if $\beta,(b_1^2 + b_2^2) \neq 0$. Applying the change of basis $P = (b_1^2 + b_2^2)^{-1} \begin{pmatrix}
b_2 & -b_1 \\
b_1 & b_2
\end{pmatrix}$, the study reduces to the stabilization of the system:
\begin{equation*}
    \begin{pmatrix}
        \dot{x}^1 \\
        \dot{x}^2
    \end{pmatrix} = \begin{pmatrix}
        \alpha & \beta \\
        - \beta & \alpha
    \end{pmatrix} \begin{pmatrix}
        x^1 \\
        x^2
    \end{pmatrix} + u\begin{pmatrix}
        0\\
        1
    \end{pmatrix} \qquad \alpha > 0 \qquad \beta \neq 0 \qquad |u| \leq 1.
\end{equation*}
We use the complex notation $z(t) = x^1(t) + ix^2(t)$. If $u = 0$ then the trajectory is a logarithmic spiral of equation $z(t)=z(0)e^{(\alpha - i \beta)t}$. If $u \pm 1$ the trajectory is a logarithmic spiral of equation $z(t)=(z(0) \pm \bar{z})e^{(\alpha - i\beta)t} \mp \bar{z}$ with $\bar{z} = \frac{\beta - i \alpha}{\alpha^2 + \beta^2}$. \\ \\
The domain of initial conditions that can be stabilized to the origin with a bounded control $|u|\leq1$ is delimited by the parameterized curves $z(t) = \pm z_{lim}e^{(\alpha - i\beta)t} \pm \bar{z},-\frac{\pi}{\beta}\leq t\leq0$ with $z_{lim} = \bar{z}(1+\frac{2}{e^{\frac{\alpha}{\beta}\pi}-1})$. A simple way to find this domain is to compute time-optimal synthesis. It should be well known and we leave the details to the reader.

\begin{prop}
    Let $x_0$ be an initial condition in the attainable set. There exists a finite time $T$ for which the lowest $L^1$-norm cost $\mu_{\infty}(x_0)$ can be reached. The corresponding optimal control $u$ and optimal trajectory $x(\cdot)$ are unique until $0_{\mathbf{R}^2}$ is reached, and the structure of $u$ is given by a switching function of the form $f(t) = r_0e^{-\alpha t }\sin(\theta_0-\beta t),r_0 > 0, \theta_0 \in [0;2\pi)$. If $|f(t)| \geq 1$, then $u(t)=sign(f(t))$, otherwise $u=0$.
\end{prop}

\noindent{\bf Proof.}
For any nontrivial solution of the adjoint equation $\dot{p}=-pA$, we obtain $p(t)B=p_2(t)=f(t)$, where $f$ is a function of the form described in the statement of the proposition. We see that any switching function has such a form. Note, that this class of functions is invariant with respect to the translation of the argument: if $f$ is in this class then
the function $t\mapsto f(t+s)$ also is, $\forall\,s\in\R$. Moreover, for such a function $f$, there is a unique $s_f\in\R$ such that $|f(s_f)| = 1$ and $|f(t)| \leq 1, \forall \: t \geq s_f$. Given $T>0$, let
$$
u(t)=\begin{cases} sign(f(t+s_f-T))&\text{if $|f(t+s_f-T)|>1$},\\
0&\text{if $|f(t+s_f-T)|\le 1$},\end{cases}
$$
and $x(t)$ be the solution of the Cauchy problem
$$
\dot x=Ax+Bu(t),\qquad x(T)=0.
$$
According to Theorem 4, $x(\cdot)$ is optimal solution for the free time problem with the initial condition $x_0=x(0)$. Moreover, any optimal solution has such a form. $\square$
\\ \\
In order to obtain images of the optimal synthesis, we use numerical simulation. We compute the trajectory given by some switching functions in reverse time and starting from the origin. The co-vector follows a trajectory of polar equation $r(\theta)=r_0e^{-\frac{\alpha}{\beta}\theta_0}(e^{\frac{\alpha}{\beta}})^{\theta}$. The numerical simulations can reduce to a one parameter simulation. For each spiral of equation $r(\theta) = C_0 (e^{\frac{\alpha}{\beta}t})^\theta$, starting from the origin in reverse time and following the spiral curve and using the ordinate of the point as the switching function, we obtain one optimal trajectory in reverse time starting from the origin. It it sufficient to make variations of $C_0$ in the interval $(0,e^{\frac{\alpha}{\beta}2\pi}]$ to obtain all possible trajectories. In red the part of the trajectory when $u=0$, in blue when $u=\pm1$.

\begin{center}
\includegraphics[scale=0.50]{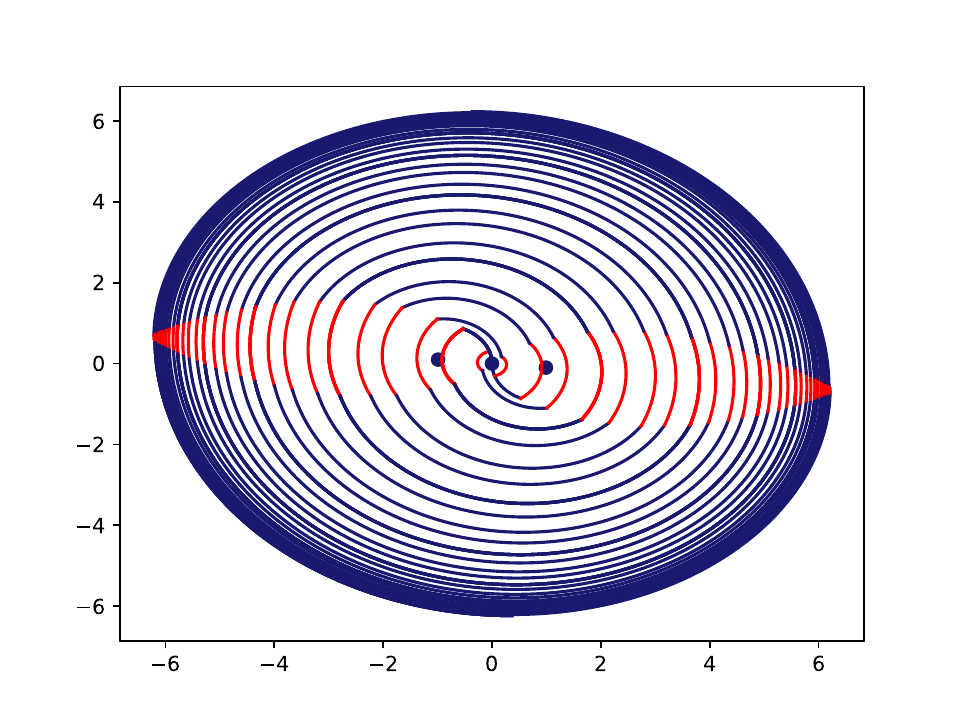}
\captionof{figure}{Some optimal trajectories for $\alpha=0.1,\beta=1$}
\end{center}

\begin{center}
\includegraphics[scale=0.50]{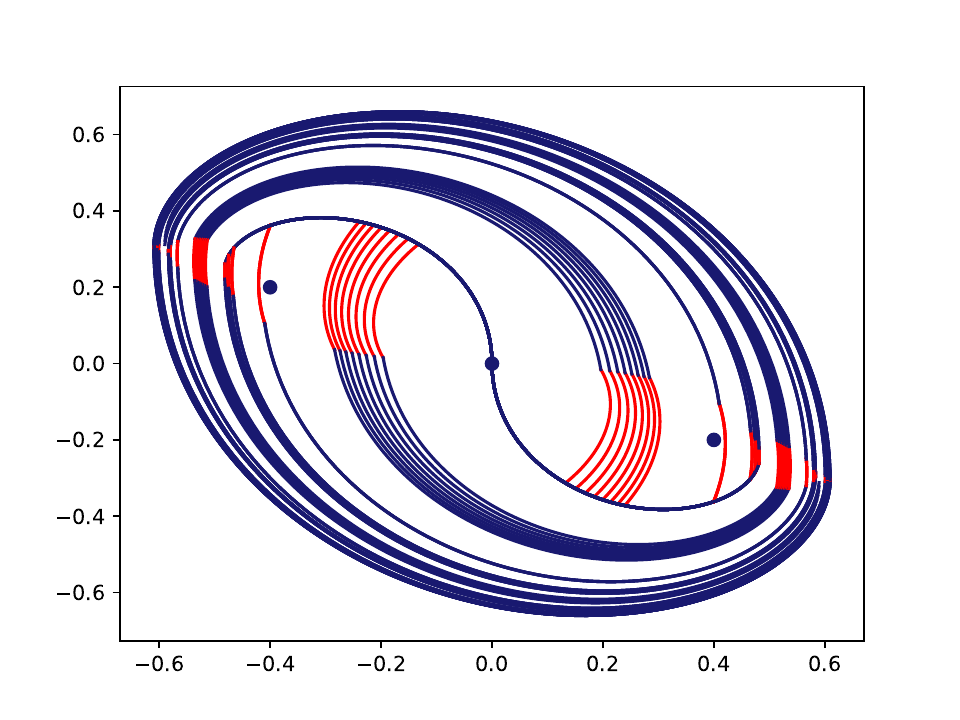}
\captionof{figure}{Some optimal trajectories for $\alpha=1,\beta=2$}
\end{center}

\section{Conclusion}

In this paper, we study optimal control problems for finite dimensional linear systems, where control parameters take values in the Euclidean unit ball and the cost is the $L^1$-norm of the control function. We study both finite and infinite horizon problems. Main results concern the existence, uniqueness and structure of optimal solutions, as well as continuous dependence of the optimal cost on the data.

\section*{Appendix}

\subsection*{A Harmonic Oscillator}

\subsection*{Proof of Theorem~5}

We consider the optimal control with the structure described in Proposition 5, we denote $\alpha_0$ the length of the last interval where $u \neq 0$ and we look at the specific case where $\delta$, the length of the time interval where $u$ is equal to $+1$ or $-1$, is equal to $\alpha_0$. \\
If on the last time interval $u \equiv 1$, then because the trajectory arrives at $0_{\mathbf{R}^2}$, the last arc of the optimal curve is an arc of the circle centered at $(1,0)$ with radius $1$, traveled clockwise. \\
After the interval $[T-\alpha_0,T]$ where $u\equiv 1$, $u \equiv 0$ for a time $\pi - \delta$ and the trajectory is a circle centered at $(0,0)$. Because here $\delta = \alpha_0$, applying the inscribed angle theorem, we obtain the upper part of the circle as the previous switching curve.
\begin{center}
\includegraphics[scale=0.25]{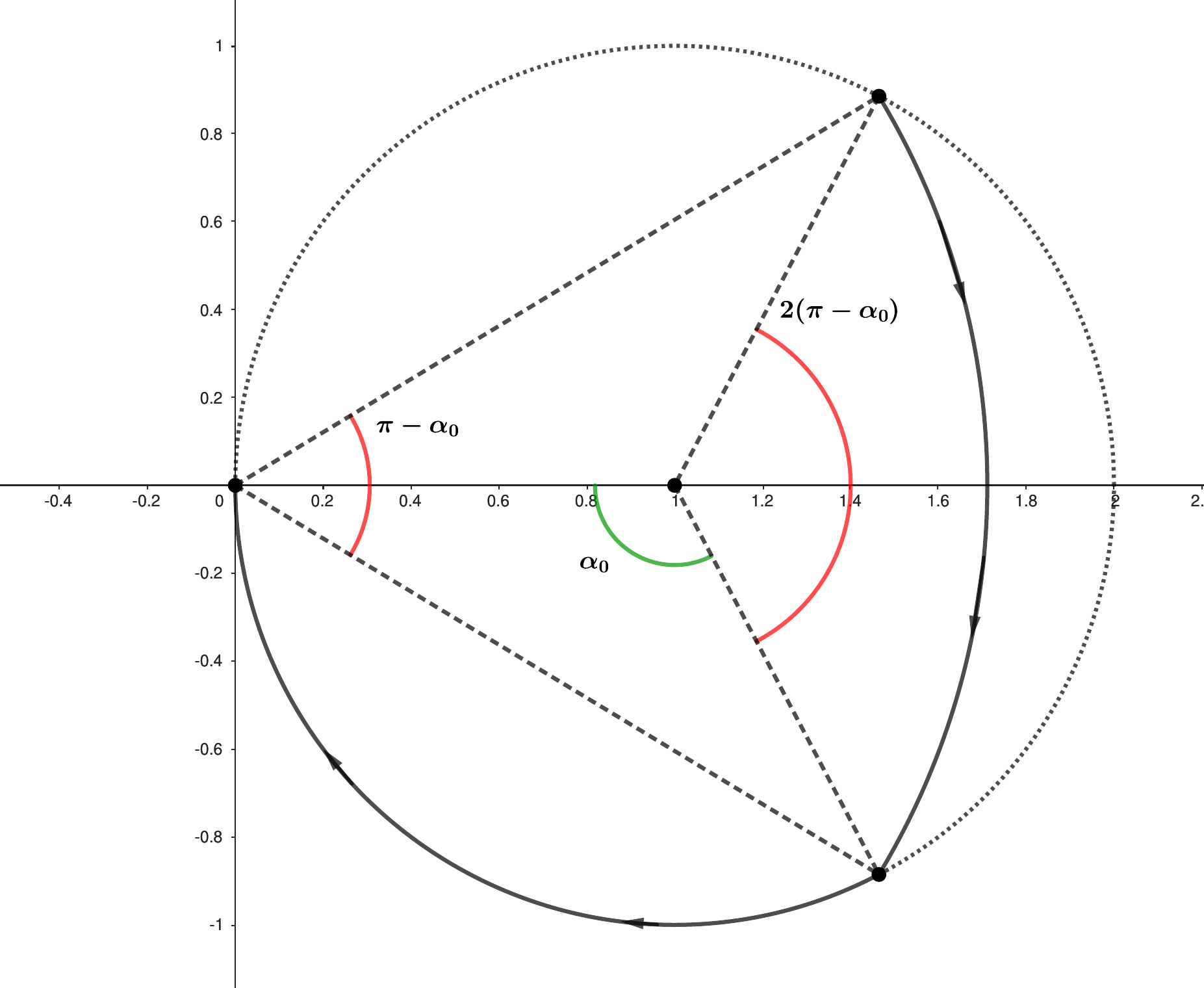}
\captionof{figure}{Last switching curve when $\delta = \alpha_0$}
\end{center}
In order to construct the last to last switching curve, we use the following geometric construction:
\begin{center}
\includegraphics[scale=0.25]{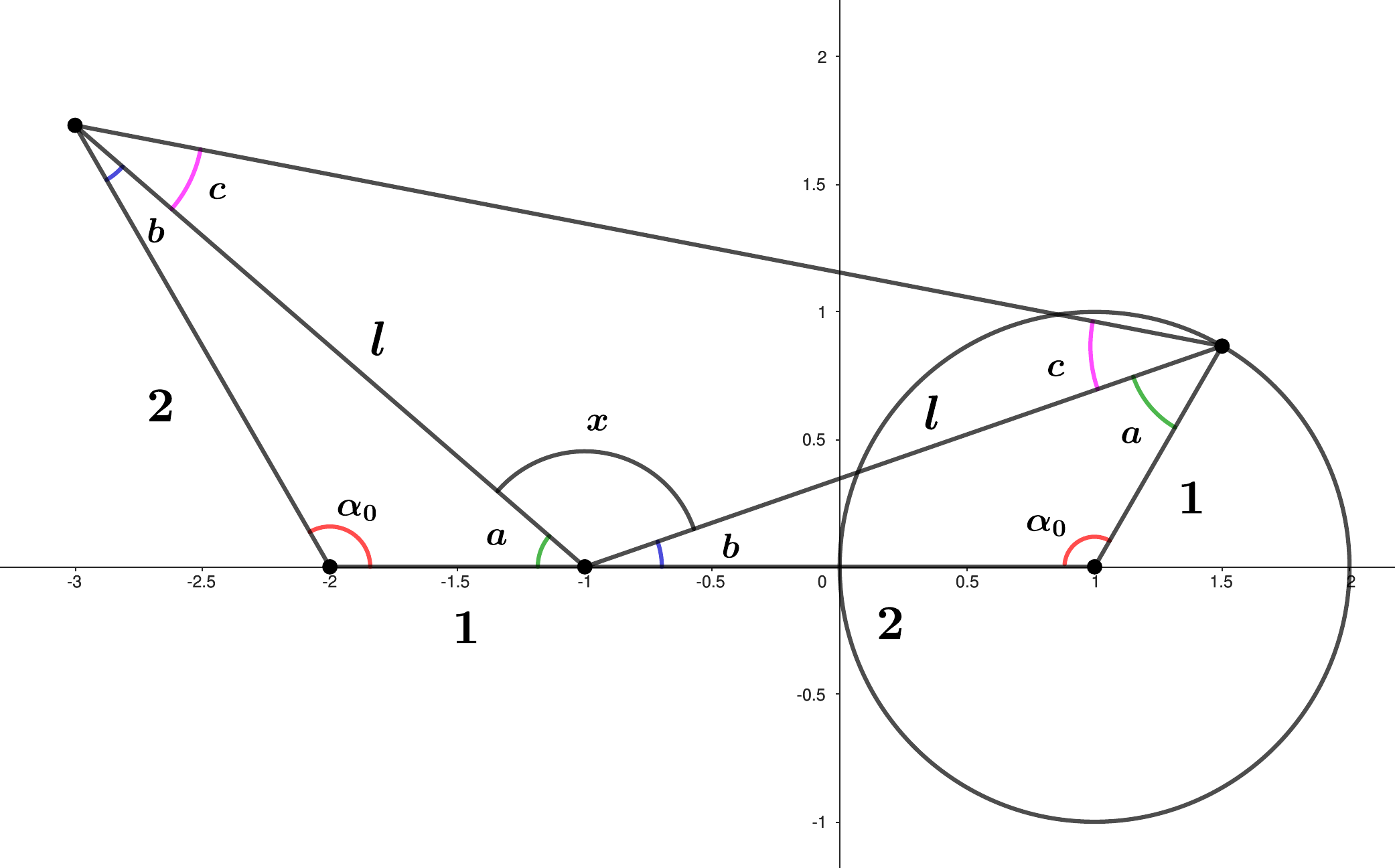}
\captionof{figure}{Construction of the last to last switching point}
\end{center}
The two triangles that have one side on the abscissa axis are the same triangles but with a different side on this axis. And so the triangle in the middle is isosceles. Then because the sum of the angles of a triangle must be equal $\pi$ we have:
\begin{equation*}
    \begin{cases}
        a + b + \alpha_0 = \pi \\
        x + 2c = \pi
    \end{cases}.
\end{equation*}
Because the sum of the angles of a quadrilateral is equal to $2\pi$ we have:
\begin{equation*}
    2\alpha_0 + a + b + 2c = 2\pi.
\end{equation*}
Finally we also note that:
\begin{equation*}
    a+b+x = \pi.
\end{equation*}
Using these equations we obtain $x=\alpha_0$. So this figure corresponds to the construction of the next to next to last switching point. The corresponding switching curve is the upper half circle centered at $(-2,0)$ with radius $2$. Repeating this scheme of construction, in the case where $u\equiv 1$ on the last time interval, we obtain all the circles centered at $(-2n,0)$ with radius $2n$ and all those centered at $(2n+1,0)$ with radius $2n+1$, where $n \in \mathbf{N}$. By symmetry, if we add the possibility of the value $u=-1$ in the last interval $[T-\alpha_0;T]$ we obtain all the circles centered at $(k,0)$ with radius $|k|$, where $k \in \mathbf{Z}$.
\begin{center}
\includegraphics[scale=0.30]{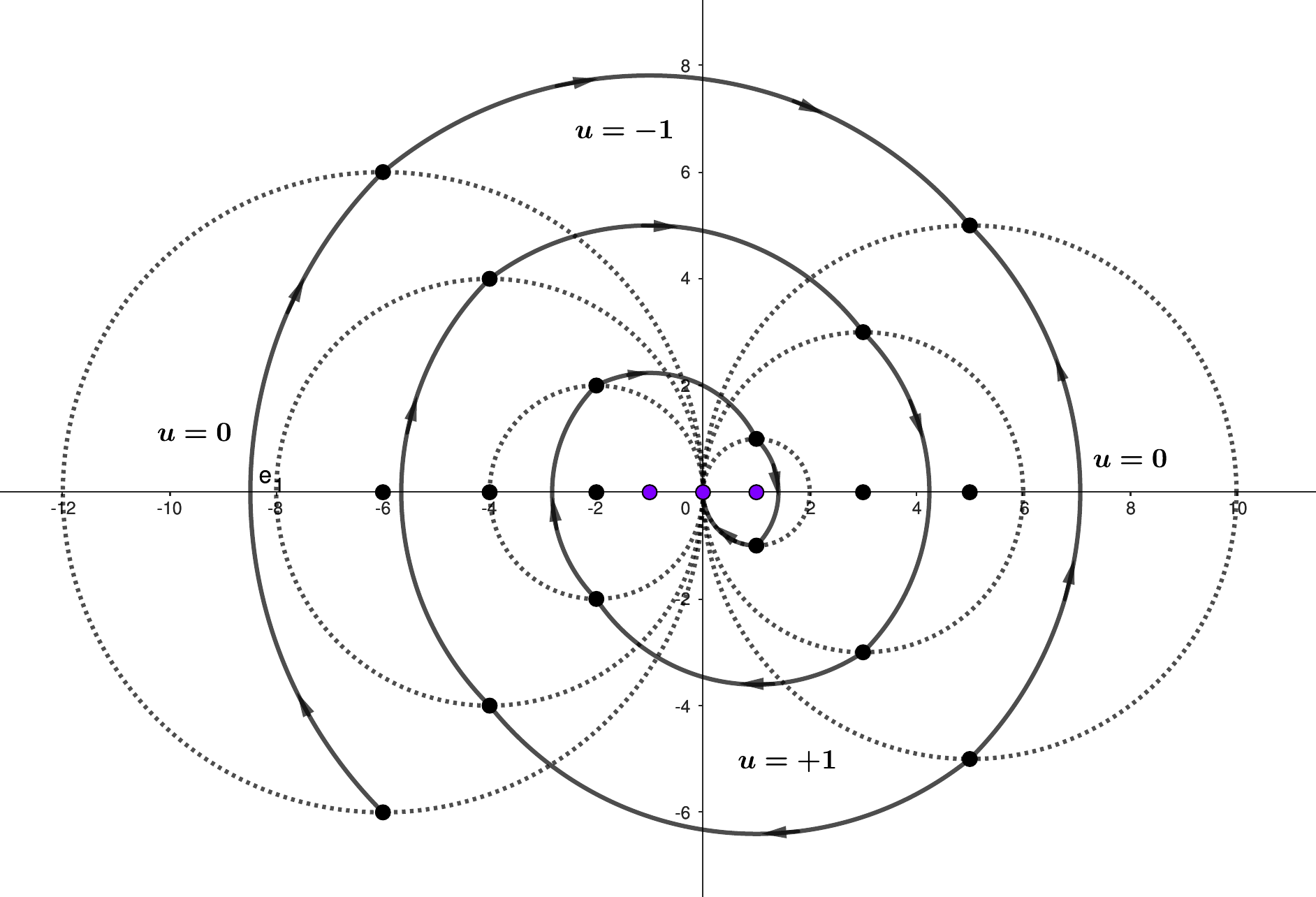}
\captionof{figure}{Switching curves when $u\equiv+1$ at the end}
\end{center}
Let $x_0 \in \mathbf{R}^2$, $x^1_{0} > 0$, $k \in \mathbf{N}, k > |x_0|$. We begin with a control $u \equiv 0$, the trajectory is a circle centered at $(0,0)$ until the first switch on the circle centered at $(k,0)$ with radius $k$. Then the trajectory is as described previously. Such a trajectory has $k$ arcs where $u\equiv+1$ or $u\equiv-1$, all of length $\delta=\alpha_0(k)$. So the $L^1$-norm cost of such a trajectory is equal to $c(k)=k\alpha_0(k)$. With the notations of the following figure, we compute $\cos(\alpha_0(k))=1-\frac{|x_0|^2}{2k^2}$.
\begin{center}
\includegraphics[scale=0.30]{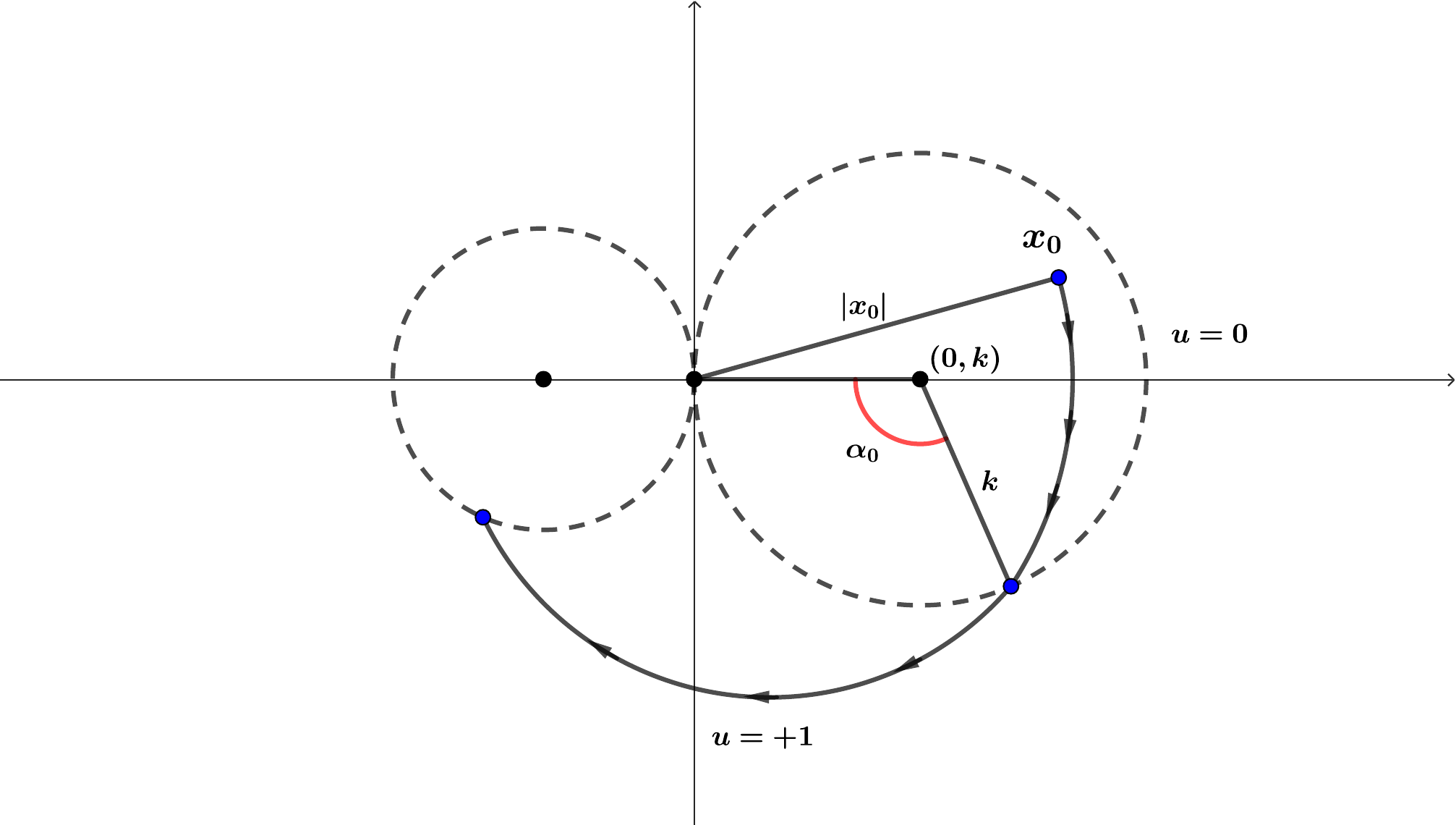}
\captionof{figure}{Construction of an optimal trajectory when time horizon increases}
\end{center}
The optimal cost $c(k)=k\arccos(1-\frac{|x_0|^2}{2k^2})$ goes to $|x_0|$ when $k$ goes to infinity. The functional $\mu_{\infty}$ is decreasing and bounded from below (by 0) so the limit exists : $\mu_{\infty}(x_0)= \lim_{T_k \rightarrow \infty} \mu_{T_k}(x_0) = |x_0|$. \\
When $k$ goes to infinity, the trajectory has more and more switches, but with arcs of circles centered at $(-1,0)$ and $(1,0)$ that are smaller and smaller.
\begin{center}
\includegraphics[scale=0.30]{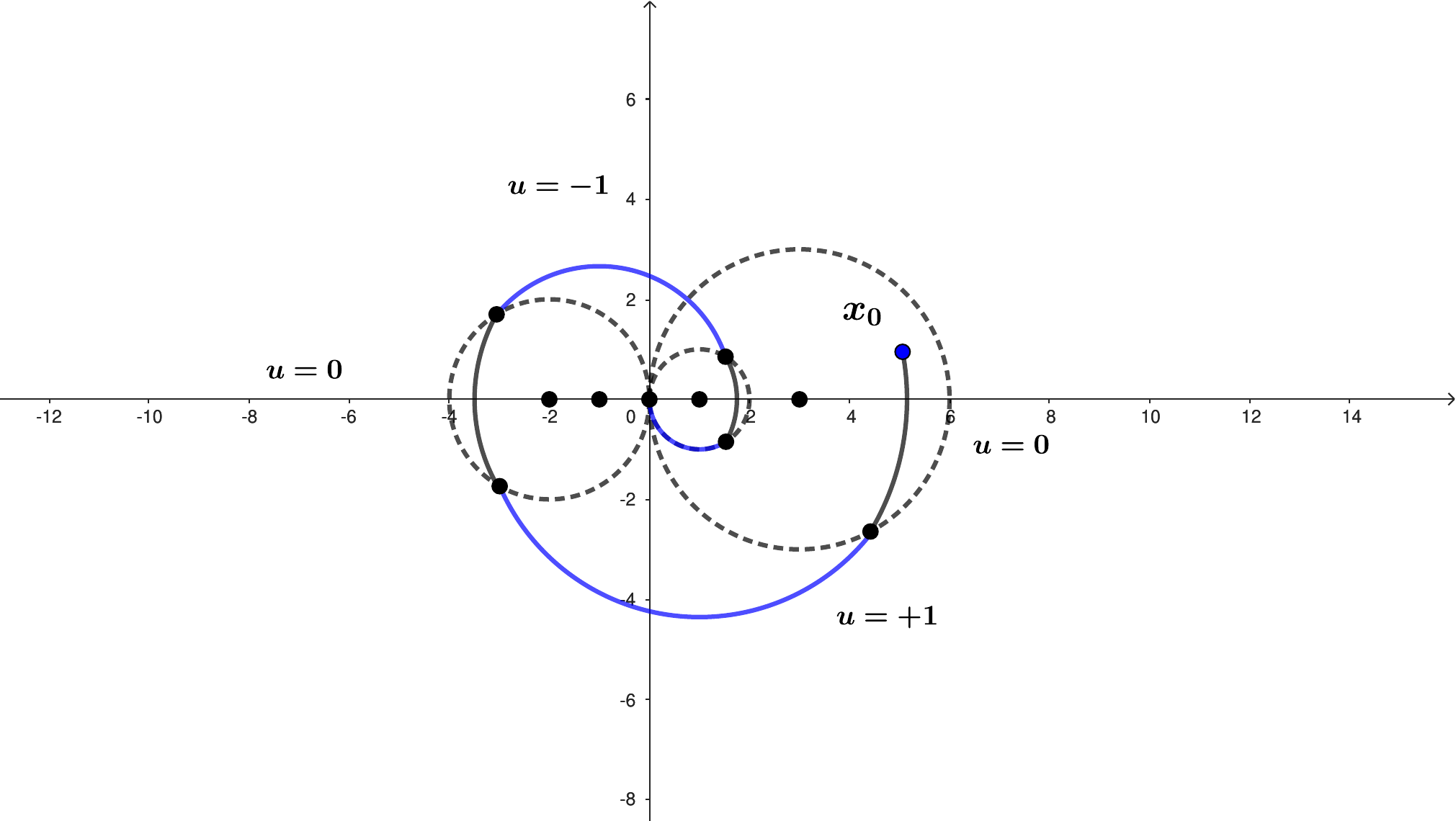}
\captionof{figure}{Trajectory for $2\pi \leq T \leq 3\pi$}
\end{center}
\begin{center}
\includegraphics[scale=0.30]{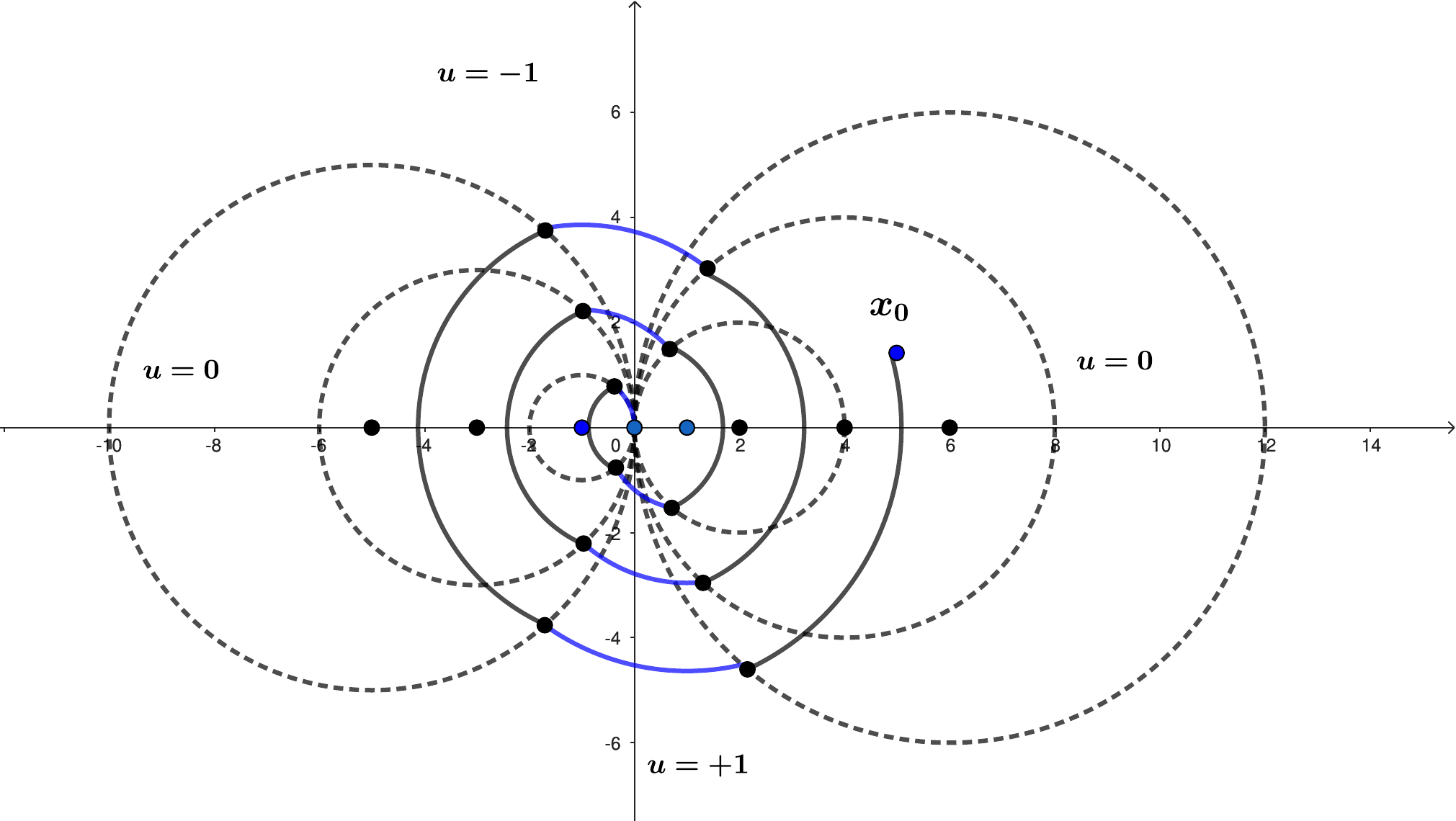}
\captionof{figure}{Trajectory for $5\pi \leq T \leq 6\pi$}
\end{center}

\end{document}